\documentclass[twoside,a4paper,10pt]{article}
\usepackage{latexsym,amssymb}
\usepackage{graphicx}
\usepackage[frenchb,english]{babel}
\usepackage[matrix,arrow]{xy}
\def \be{\begin{eqnarray*}}
\def \ee{\end{eqnarray*}}
\def \ben{\begin{enumerate}}
\def \een{\end{enumerate}}
\def \beit{\begin{itemize}}
\def \eeit{\end{itemize}}
\def \bui#1#2{\mathrel{\mathop{\kern 0pt#1}\limits^{#2}}}
\def \buil#1#2{\mathrel{\mathop{\kern 0pt#1}\limits_{#2}}}
\def \bfll{\begin{flushleft}}
\def \efll{\end{flushleft}}
\def \bflr{\begin{flushright}}
\def \eflr{\end{flushright}}

\def \lra{\longrightarrow}
\def \lmt{\longmapsto}
\def \ovl{\overline}

\def \wit{\widetilde}
\def \wnabla{\wit{\nabla}}

\def \Id{\mathrm{Id}}

\def \C{\mathbb{C}}
\def \SO{\mathrm{SO}}

\def \s{\mathbb{S}}
\def \SU{\mathrm{SU}}

\def \R{\mathbb{R}}

\newcommand{\pa}[1]{\left(#1\right)}

\newtheorem{ethm}{Theorem}[section]
\newtheorem{edefi}[ethm]{Def\mbox{}inition}
\newtheorem{elemme}[ethm]{Lemma}
\newtheorem{erem}[ethm]{Note}
\newtheorem{erems}[ethm]{Notes}
\newtheorem{ecor}[ethm]{Corollary}
\newtheorem{prop}[ethm]{Proposition}

\newtheorem{eexemples}[ethm]{Examples}
\title{Geometric aspects of transversal   
Killing spinors on Riemannian f\mbox{}lows}
\author{Nicolas Ginoux\footnote{Institut f\"ur Mathematik - Geometrie,
Universit\"at Potsdam,
Am Neuen Palais 10,
D-14469 Potsdam,
E-mail: ginoux@math.uni-potsdam.de}, Georges Habib\footnote{Max-Planck Institut f\"ur Mathematik in den Naturwissenschaften, Inselstrasse 22, D-04103 Leipzig, Germany, E-mail: habib@mis.mpg.de}}

\begin{document}
\date{}
\maketitle

\noindent\begin{center}\begin{tabular}{p{115mm}}
\begin{small}{\bf Abstract.} We study a Killing spinor type equation on spin Riemannian f\mbox{}lows. We prove integrability conditions and partially classify those Riemannian f\mbox{}lows $M$ carrying non-trivial solutions to that equation in case $M$ is a local Riemannian product, a Sasakian manifold or $3$-dimensional.
\end{small}\\
\end{tabular}\end{center}
 $ $\\

\noindent\begin{small}{\it Mathematics Subject Classif\mbox{}ication}: 53C12, 53C27 \\
\end{small}

\noindent\begin{small}{\it Key words}: Foliations, Spin Geometry\\
\end{small}

\section{Introduction}
\noindent Killing spinors on Riemannian spin manifolds are smooth sections of the spinor bundle of which covariant derivative is proportional to the Clif\mbox{}ford multiplication. Those manifolds carrying non-zero Killing spinors have been well understood for a long time, see e.g. \cite{BHMM} for a survey. Such a manifold is Einstein; if it is furthermore compact and non Ricci f\mbox{}lat, then it cannot carry any non-zero non-trivial parallel form (see e.g. \cite{Hij98}). In particular, the existence of non-zero Killing spinors imposes very rigid conditions to the geometry of the underlying manifold.\\ 

\noindent In this article we transpose the Killing spinor equation to the set up of Riemannian f\mbox{}lows \cite{Car}, which roughly speaking, are local Riemannian submersions with $1$-dimensional fibres (see Section \ref{sprelim} for the definition). These contain among others all $\s^1$-bundles with totally geodesic f\mbox{}ibres, all manifolds carrying a non-zero parallel form and all Sasakian manifolds \cite{BGal99}. In the def\mbox{}inition of that equation - that we call transversal Killing spinor - we allow the f\mbox{}irst derivatives of the spinor f\mbox{}ield to behave dif\mbox{}ferently along the leaves and along the orthogonal distribution of the foliation respectively, see Def\mbox{}inition \ref{ddefKT}.  There are several motivations for this study. Historically those spinor f\mbox{}ields f\mbox{}irst appear in the limiting case of an eigenvalue estimate proved by B. Alexandrov, G. Grantcharov and S. Ivanov for the Dirac operator on compact Riemannian spin manifolds with a parallel $1$-form \cite[eq. (8)]{AlexGrantchIv98}. On the other hand, they stand for the most natural tools in the study of the spectrum of the Dirac operator on submersions over real space forms \cite{Habibthese,GinHab2}.\\

\noindent The paper is organized as follows. In the f\mbox{}irst part we recall basic facts on spin Riemannian f\mbox{}lows. In the second one, we prove integrability conditions for spin Riemannian f\mbox{}lows admitting non-zero transversal Killing spinors, see Theorem \ref{tcnexistpsi}. It should be noticed that the resulting geometric conditions hold up to homothetic deformations of the metric along the leaves (see Lemma \ref{lexistspinKTdeform}). We describe then in Proposition \ref{pS1fibres} important examples of Riemannian f\mbox{}lows carrying transversal Killing spinors, which arise as submersions over manifolds with Killing spinors. E\-xa\-mi\-ning group-equivariance conditions we formulate in Section \ref{exakillspi} more precise statements in case the f\mbox{}low is a local Riemannian product or f\mbox{}lat and which we illustrate on $3$-dimensional Bieberbach manifolds \cite{Pfaeffle00}. In Section \ref{sasaki}, we translate the results of Theorem \ref{tcnexistpsi} in the Sasakian setting. In that case transversal Killing spinors can be related to classical ones, see Proposition \ref{pSasakiKilling}. As an example, we describe all transversal Killing spinors on the Berger spheres (here beware of our def\mbox{}inition of Berger spheres, see Subsection \ref{subsBerger}). In the last section, we restrict ourselves to $3$-dimensional Riemannian f\mbox{}lows and simplify the conditions of Theo\-rem \ref{tcnexistpsi}. We end with the classif\mbox{}ication of compact $3$-dimensional $\eta$-Einstein minimal Riemannian f\mbox{}lows carrying non-zero transversal Killing spinors.\\

\noindent {\bf Acknowledgement.} The authors would like to thank the universities of Nancy and Potsdam as well as the Sonderforschungsbereich 647 ``{\it Raum - Zeit - Materie. Analytische und Geometrische Strukturen}'' of the Deutsche Forschungsgemeinschaft for their support in the preparation of \cite{GinHab2} and this paper. It's also a pleasure to thank Bernd Ammann, Christian B\"ar and Oussama Hijazi for valuable comments.

\section{Preliminaries}\label{sprelim}

\noindent For preliminaries on Riemannian spin foliations we refer to \cite[Chap. 1]{Habibthese}.\\

\noindent Throughout the paper the triple $(M^{n+1},g,\mathcal{F})$ will denote an $(n+1)$-dimensional Riemannian manifold endowed with a Riemannian f\mbox{}low $\mathcal{F}$ given by a smooth unit vector f\mbox{}ield $\xi$ \cite{Car}. That means $\mathcal{F}$ defines a $1$-dimensional foliation of $M$ satisfying for all $Z,W$ orthogonal to $\xi$ the relation $(\mathcal{L}_\xi g)(Z,W)=0$ \cite{Rei} where $\mathcal{L}_\xi$ is the Lie derivative in the direction of $\xi$. Recall from the hypothesis on the f\mbox{}low being Riemannian, the endomorphism f\mbox{}ield $h:=\nabla^M\xi$ (known as the O'Neill tensor \cite{One}) of the normal bundle $Q=\xi^\perp$ is skew-symmetric w.r.t. the induced metric $g$.
We denote by $\kappa:=\nabla_\xi^M\xi$ the mean curvature of the f\mbox{}low.  The f\mbox{}low is called minimal if $\kappa=0$, which is equivalent to the integral curves of $\xi$ being geodesics or to the stronger assumption $\xi$ being a Killing vector f\mbox{}ield on $(M,g)$. Since $\xi$ has length $1$, we have $\kappa\in\Gamma(Q)$ and $g(h(Z),\xi)=0$ for every $Z \in Q$ (hence $h$ actually maps $Q$ into $Q$). In particular one may associate a $2$-form $\Omega$ to $h$ on $Q$ through $\Omega(Z,W):=g(h(Z),W)$ for all sections $Z,W\in \Gamma(Q)$.\\

\noindent Moreover, the normal bundle $Q\rightarrow M$ carries a natural covariant derivative $\nabla$ def\mbox{}ined for every section $Z$ of $Q$ by \cite{Ton}
\[\nabla_XZ:=\left|\begin{array}{ll}[\xi,Z]^Q&\textrm{ if }X=\xi\\\\
(\nabla_X^MZ)^Q&\textrm{ if }X\in \Gamma(Q),\end{array}\right. \]
where $\nabla^M$ denotes the Levi-Civita connection on $TM$ and $(\cdot)^Q$ the orthogonal projection onto $Q\subset TM$. From its def\mbox{}inition the covariant derivative $\nabla$ can be expressed in terms of the Levi-Civita connection $\nabla^M$ through the following Gauss-type formula

\begin{equation}\label{eqGaussQ}
\left|\begin{array}{ll}
\nabla_\xi^MZ&= \nabla_\xi Z+h(Z)-g(Z,\kappa)\xi\\\\
\nabla_Z^MW&= \nabla_ZW-g(h(Z),W)\xi
\end{array}\right.
\end{equation}
for every $Z,W\in\Gamma(Q)$.\\

\noindent Since $TM=\R\xi\oplus\, Q$ the normal bundle $Q\rightarrow M$ is spin (as a vector bundle) if and only if $TM$ is, and in that case spin structures on $TM$ and $Q$ respectively are in one-to-one correspondence. If we assume $M$ to be spin and carry a f\mbox{}ixed spin structure, then so does $Q$. In that case $M$ carries its own spinor bundle $\Sigma M\rightarrow M$ as well as the spinor bundle of $Q$, that we denote by $\Sigma Q\rightarrow M$. Actually there exists a unitary isomorphism of Hermitian vector bundles (that we denote by the identity map $\varphi\mapsto \varphi$)
\[\Sigma M\lra\left|\begin{array}{ll}\Sigma Q&\textrm{ if }n\textrm{ is even}\\
 \Sigma Q\oplus \Sigma Q&\textrm{ if }n\textrm{ is odd}\end{array}\right.\]
satisfying, for every $Z\in \Gamma(Q)$ and $\varphi\in\Gamma(\Sigma M)$ \cite{Baer98, Mor01}
\beit\item W.r.t. the Clif\mbox{}ford multiplications ``$\cdot$'' in $\Sigma M$ and ``$\buil{\cdot}{Q}$'' in $\Sigma Q$ respectively 
\[\xi\cdot Z\cdot\varphi=\left|\begin{array}{ll} Z\buil{\cdot}{Q}\varphi&\textrm{ if }n\textrm{ is even}\\
(Z\buil{\cdot}{Q}\oplus-Z\buil{\cdot}{Q})\varphi&\textrm{ if }n\textrm{ is odd.}\end{array}\right.\]
\item The Clif\mbox{}ford action of $i\xi$ is given by
\[i\xi\cdot=\left|\begin{array}{ll}\mathrm{Id}_{\Sigma^+Q}\oplus-\mathrm{Id}_{\Sigma^-Q}&\textrm{if }n\textrm{ is even}\\
&\\
\left(\begin{array}{ll}0&\mathrm{Id}_{\Sigma Q}\\\mathrm{Id}_{\Sigma Q}&0\end{array}\right)&\textrm{if }n\textrm{ is odd.}\end{array}\right.\]
\item W.r.t. the spinorial Levi-Civita connections $\nabla^M$ on $\Sigma M$ and $\nabla$ on $\Sigma Q$ one has \cite[eq. (2.4.7)]{Habibthese}
\begin{equation}\label{eqrestrictspin}
\left|\begin{array}{ll}\nabla_\xi^M\varphi&= \nabla_\xi\varphi+\frac{1}{2}\Omega\cdot\varphi+\frac{1}{2}\xi\cdot\kappa\cdot\varphi\\
 &\\
\nabla_Z^M\varphi&= \nabla_Z\varphi+\frac{1}{2}\xi\cdot h(Z)\cdot\varphi.\end{array}\right.
\end{equation}

\eeit
In particular the covariant derivative $\nabla$ is metric on $\Sigma M$: if $\langle\cdot\,,\cdot\rangle$ denotes a natural Hermitian inner product on $\Sigma M$, then $X\langle\varphi,\psi\rangle=\langle\nabla_X\varphi,\psi\rangle+\langle\varphi,\nabla_X\psi\rangle$ for all $X\in \Gamma(TM)$ and $\varphi,\psi\in\Gamma(\Sigma M)$. It also follows from (\ref{eqrestrictspin}) that the Clif\mbox{}ford action of $\xi$ is $\nabla$-parallel: $\nabla_X(\xi\cdot\varphi)=\xi\cdot\nabla_X\varphi$ for every $X\in \Gamma(TM)$ and $\varphi\in\Gamma(\Sigma M)$. Therefore, if one def\mbox{}ines $\Sigma_+M$ and $\Sigma_-M$ by
\[\Sigma_\pm M:=\mathrm{Ker}\Big(i\xi\cdot\mp\mathrm{Id}_{\Sigma M}\Big),\]
then $\Sigma M$ splits into the orthogonal and $\nabla$-parallel direct sum $\Sigma M=\Sigma_+M\oplus \Sigma_-M$. Furthermore both $\Sigma_+M$ and $\Sigma_-M$ have the same rank since they are exchanged by the Clif\mbox{}ford action of any non-zero section $Z\in \Gamma(Q)$. In the case where $n$ is even, one has $\Sigma_\pm M=\Sigma^\pm Q$, however in the case where $n$ is odd $\Sigma_\pm M$ never coincides with one of the two copies of $\Sigma Q$ hence with one of the eigenspaces $\Sigma^\pm M$ of the Clif\mbox{}ford action of the complex volume form of $M$.\\

\noindent In the following almost all Riemannian f\mbox{}lows under consideration will be mi\-ni\-mal, i.e., $\xi$ will be a (unit) Killing vector f\mbox{}ield. More precisely, we shall mainly deal with the following families of Riemannian f\mbox{}lows, which of course are not disjoint from each other: the case where $\xi$ is the fundamental vector f\mbox{}ield of a free isometric $\s^1$-action with totally geodesic orbits, the case where $\xi$ is parallel, corresponding to local Riemannian products of a one-dimensional manifold with an $n$-dimensional one (this is also equivalent to $\kappa=0$ and $h=0$), and the case where $\xi$ is the Reeb vector f\mbox{}ield of a Sasakian manifold 

\begin{edefi}\label{defSasaki}
A Riemannian manifold is called \emph{Sasakian} if and only if it is a Riemannian f\mbox{}low $(M^{n+1},g,\mathcal{F})$ satisfying 
\beit\item[i)] $\kappa=0$, i.e., the f\mbox{}low is minimal,
\item[ii)] $h^2=-\mathrm{Id}_Q$, i.e., $h$ is an almost-Hermitian structure on $Q$,
\item[iii)] $\nabla h=0$, i.e., $h$ is parallel on $Q$ (hence is a K\"ahler structure on $Q$).
\eeit
\end{edefi} 

\noindent It can be easily checked that this def\mbox{}inition is equivalent to the usual one, where one requires $\xi$ to be a unit Killing vector f\mbox{}ield satisfying
\[
\left|\begin{array}{ll}(\nabla^M\xi)^2&=-\mathrm{Id}_{TM}+\xi^\flat\otimes \xi\\
(\nabla_X^M\nabla^M\xi)(Y)&=g(\xi,Y)X-g(X,Y)\xi\end{array}\right.
\]
for all $X,Y\in \Gamma(TM)$. From Definition \ref{defSasaki} the normal bundle  $Q$ of any Sasakian manifold carries a canonical K\"ahler structure. In particular such a manifold is always odd-dimensional. We shall from now on denote $m:=\frac{n}{2}$. In the following we shall also omit to write $\mathcal{F}$ for the flow and consider a Sasakian manifold as a triple $(M^{2m+1},g,\xi)$.\\\\ 
If now a Sasakian manifold $M$ is spin, then the Clif\mbox{}ford action of the $2$-form $\Omega$ (which is then the K\"ahler form of $Q$) splits $\Sigma M$ into the following orthogonal and $\nabla$-parallel decomposition \cite{Kir86}:
\begin{equation}\label{eqdecSigmaSasaki}
\Sigma M=\bigoplus_{r=0}^{m}\Sigma_r M,
\end{equation}
where $\Sigma_r M$ is the eigenbundle associated with the eigenvalue $i(2r-m)$ of $\Omega$ for every $r\in\{0,\ldots,m\}$. W.r.t. the Clif\mbox{}ford action of $i\xi$ one has $\Sigma_+M=\bigoplus_{r\,\textrm{\scriptsize even}}\Sigma_rM$ and $\Sigma_-M=\bigoplus_{r\,\textrm{\scriptsize  odd}}\Sigma_rM$, that is $i\xi\cdot_{|_{\Sigma_rM}}=(-1)^r\mathrm{Id}_{\Sigma_rM}$ \cite{FK}.
Moreover, the subspaces $\Sigma_0M$ and $\Sigma_{m}M$ can be characterized by the property that for all $Z$ orthogonal to $\xi$ we have 
\begin{equation} \label{eqidentsasa}
h(Z)\cdot\psi_0=iZ\cdot\psi_0  \quad{\rm and}\quad h(Z)\cdot\psi_m=-iZ\cdot\psi_m.
\end{equation}

\noindent Replacing the metric by a positive scalar multiple of it obviously preserves the structure of Riemannian flow. There exists however a less trivial type of flow-preserving deformations of the metric that we will need in the next sections and which are called $\mathcal{D}$-homothetic:

\begin{edefi}[S. Tanno \cite{Tan68}]\label{ddefDhomot}
Let $(M^{n+1},g,\mathcal{F})$ be a Riemannian f\mbox{}low. given by a unit vector field $\xi$. A \emph{$\mathcal{D}$-homothetic deformation} of $g$ is a metric $\ovl{g}_t$ on $M$ of the form
\[\ovl{g}_t:=t^2g_{|_{\R\xi}}+tg_{|_{Q}}\]
for some real number $t>0$.
\end{edefi}

\noindent A $\mathcal{D}$-homothetic deformation of $g$ may be obtained as follows: f\mbox{}irst rescale $g$ by a factor $t$ in the direction of the f\mbox{}low, then multiply the obtained metric by $t$. It is f\mbox{}irst to be noticed that a $\mathcal{D}$-homothetic deformation of a Riemannian f\mbox{}low is again a Riemannian f\mbox{}low, more precisely:

\begin{elemme}\label{ldeformDhomot}
Let $(M^{n+1},g,\mathcal{F})$ be a Riemannian f\mbox{}low given by a unit vector field $\xi$ and $\ovl{g}_t:=t^2g_{|_{\R\xi}}+tg_{|_{Q}}$ be a $\mathcal{D}$-homothetic deformation of $g$ with $t\in\R_+^*$. Then $(M^{n+1},\ovl{g}_t,\mathcal{F})$ is a Riemannian f\mbox{}low described by the unit vector field $\ovl{\xi}^t:=\frac{1}{t}\xi$. If furthermore $\ovl{\xi}^t$, $\ovl{\nabla}^t$, $\ovl{h}^t$, $\ovl{\kappa}^t$ denote the corresponding objects for $\ovl{g}_t$, then the following holds:
\beit
\item[i)] One has $\ovl{\xi}^t=\frac{1}{t}\xi$, $\ovl{h}^t=h$ and $\ovl{\kappa}^t=\frac{1}{t}\kappa$.
\item[ii)] On $Q$ one has $\ovl{\nabla}^t=\nabla$.
\item[iii)] If furthermore $M$ is spin, then there exists a unitary isomorphism
\be
\Sigma_gM&\lra&\Sigma_{\ovl{g}_t}M\\
\varphi&\lmt&\ovl{\varphi}^t
\ee
s.t., if ``\,$\ovl{\cdot}$'' denotes the Clif\mbox{}ford multiplication on $\Sigma_{\ovl{g}_t}M$,
\beit\item[$\bullet$] $\ovl{\xi\cdot\varphi}=\ovl{\xi}^t\,\ovl{\cdot}\,\ovl{\varphi}^t$ and $\ovl{Z\cdot\varphi}=\frac{1}{\sqrt{t}}Z\,\ovl{\cdot}\,\ovl{\varphi}^t$ for every $Z\in \Gamma(Q)$.
\item[$\bullet$] $\ovl{\nabla}^t_X\ovl{\varphi}^t=\ovl{\nabla_X\varphi}^t$ for every $X\in \Gamma(TM)$.
\eeit
\item[iv)] If $(M^{n+1},g,\xi)$ is Sasakian, then so is $(M^{n+1},\ovl{g}_t,\ovl{\xi}_t)$.
\eeit
\end{elemme}

\noindent The proof of Lemma \ref{ldeformDhomot} consists of elementary computations and identif\mbox{}ications that we leave to the reader.

\section{Transversal Killing spinors}
\subsection{Def\mbox{}inition}
\noindent We generalize in some sense the Killing spinor equation (see e.g. \cite{BFGK,Hij98} for re\-fe\-ren\-ces on that topic) to the set up of Riemannian f\mbox{}lows.

\begin{edefi}\label{ddefKT}
Let $(M^{n+1},g,\mathcal{F})$ be a spin Riemannian f\mbox{}low given by a unit vector field $\xi$. Let $\alpha,\beta\in\C$. An \emph{$(\alpha,\beta)$-transversal Killing spinor} on $M$ is a smooth section $\psi$ of $\Sigma M$ satisfying, for all $Z\in \Gamma(Q)$,
\begin{equation}\label{eqkilling}
\left|\begin{array}{ll}
\nabla_{\xi}\psi&=\alpha\,\xi\cdot\psi\\
 &\\
\nabla_Z\psi&=\beta\, \xi\cdot Z\cdot\psi.
\end{array}\right.
\end{equation}
If $\alpha=0$, then $\psi$ is called a \emph{basic $\beta$-Killing spinor} (see also \cite{Habibthese}), and if $\alpha=\beta=0$ it is called \emph{basic parallel} or \emph{transversally parallel spinor}.
\end{edefi}

\noindent F\mbox{}irst note that an $(\alpha,\beta)$-transversal Killing spinor is a parallel section of $\Sigma M$ w.r.t. the covariant derivative $\wnabla$ def\mbox{}ined by
\[\wnabla_X\varphi:=\nabla_X\varphi-\alpha g(X,\xi)\xi\cdot\varphi-\beta\xi\cdot X\cdot\varphi-\beta g(X,\xi)\varphi\]
for all $X\in \Gamma(TM)$ and $\varphi\in\Gamma(\Sigma M)$.
Hence if an $(\alpha,\beta)$-transversal Killing spinor vanishes at one point it vanishes everywhere on $M$.

\begin{erems}\label{remdefKT}
{\rm\noindent\ben\item If $\psi$ is an $(\alpha,\beta)$-transversal Killing spinor on $(M^{n+1},g,\mathcal{F})$, then $\xi\cdot\psi$ is an $(\alpha,-\beta)$-transversal Killing spinor. Therefore $\beta$ can always be changed into $-\beta$, independently of the dimension or the orientation of the manifold. This is in general not possible for $\alpha$, see e.g. Notes \ref{remsphereBerger}.
\item Let $\psi$ be an $(\alpha,\beta)$-transversal Killing spinor and $\psi=\psi_++\psi_-$ its decomposition w.r.t. the Clif\mbox{}ford action of $i\xi$ (i.e., $i\xi\cdot\psi_\pm=\pm\psi_\pm$). Then $\psi_\pm$ satisf\mbox{}ies $\nabla_\xi\psi_\pm=\alpha\xi\cdot\psi_\pm$ and $\nabla_Z\psi_\pm=\beta\xi\cdot Z\cdot\psi_\mp$ for all $Z\in \Gamma(Q)$ (this follows from the action of $\xi$ being $\nabla$-parallel and from $\xi\cdot Z=-Z\cdot\xi$ for all $Z\in \Gamma(Q)$). Therefore $\psi_\pm$ is again a transversal Killing spinor only if $\beta=0$; if $\beta\neq 0$, then $\psi$ vanishes as soon as $\psi_+$ or $\psi_-$ vanishes on a non-empty open subset of $M$.
\item If $\alpha$ and $\beta$ are real, then any $(\alpha,\beta)$-transversal Killing spinor has constant length on $M$, since in that case it can be easily checked that the covariant derivative $\wnabla$ above is metric. 
\item In the particular case where $h=0$ and $\kappa=0$, i.e., $M$ carries a parallel unit vector f\mbox{}ield, then $(0,\beta)$-transversal Killing spinors for some $\beta\in\R$ are exactly the spinor f\mbox{}ields def\mbox{}ined by B. Alexandrov, G. Grantcharov and S. Ivanov in \cite[eq. (8)]{AlexGrantchIv98} and studied in \cite[Thm. 3.1]{AlexGrantchIv98}.
\item We notice that Th. Friedrich and E. C. Kim defined on a Sasakian manifold $M$ the notion of {\it quasi-Killing spinor of type $(a,b)$} \cite[page 23]{FK} which are $a$-Killing (for the Levi-Civita connection on $M$) in the direction of the normal bundle $Q$ and $(a+b)$-Killing in the direction of $\xi$. They show that the condition of the f\mbox{}low being $\eta$-Einstein is sufficient and necessary for the existence of such spinors. In that case and for a suitable choice of $a$ and $b$ (the product $a(a+b)$ could be not zero) they are solutions of the so-called {\it Einstein-Dirac equations}. Here we point out that the notion of transversal Killing spinors is in general different from the quasi-Killing spinors since in our consideration it is Killing for the transversal connection.  
\item A similar equation appears in \cite[page 137 eq. (8.3)]{Moroithese}, where however the connection $\nabla$ denotes the Levi-Civita connection of $\Sigma M$ (in our notation it corresponds to $\nabla^M$).
\een
}\end{erems}



\noindent In the following subsections we want to characterize those Riemannian f\mbox{}lows that admit non-trivial $(\alpha,\beta)$-transversal Killing spinors. The following lemma follows straightforward from Lemma \ref{ldeformDhomot}.

\begin{elemme}\label{lexistspinKTdeform}
Let $(M^{n+1},g,\mathcal{F})$ be a spin Riemannian f\mbox{}low. For a f\mbox{}ixed $t\in\R_+^*$ let $\ovl{g}_t:=t^2g_\xi+tg_{|_{ Q}}$ be a $\mathcal{D}$-homothetic deformation of $g$.
If $\psi$ is an $(\alpha,\beta)$-transversal Killing spinor on $(M^{n+1},g,\mathcal{F})$, then $\ovl{\psi}^t$ is an $(\ovl{\alpha}^t,\ovl{\beta}^t)$-transversal Killing spinor on $(M^{n+1},\ovl{g}_t,\ovl{\xi}^t:=\frac{1}{t}\xi)$ with $\ovl{\alpha}^t:=\frac{\alpha}{t}$ and $\ovl{\beta}^t:=\frac{\beta}{\sqrt{t}}$.
\end{elemme}

\subsection{General integrability conditions for transversal Killing spinors}

\begin{ethm}\label{tcnexistpsi}
Let $(M^{n+1},g,\mathcal{F})$ be a spin Riemannian f\mbox{}low given by a unit vector field $\xi$ carrying an $(\alpha,\beta)$-transversal Killing spinor $\psi$. Let $\mathrm{Ric}_M$ and $\mathrm{Scal}_M$ denote the Ricci tensor and the scalar curvature of $(M,g)$ respectively. Then for any local orthonormal basis $\{e_j\}_{1\leq j\leq n}$ of $Q$ one has
\be
\mathrm{Ric}_M(\xi)\cdot\psi&=&(|h|^2-|\kappa|^2)\xi\cdot\psi+4n\alpha\beta\psi+2\alpha\xi\cdot\kappa\cdot\psi+\kappa\cdot\Omega\cdot\psi\\
& &+4h(\kappa)\cdot\psi-\sum_{j=1}^n\xi\cdot e_j\cdot\nabla_{e_j}^M\kappa\cdot\psi\\
& &+\frac{1}{2}\sum_{j,k=1}^n e_j\cdot e_k\cdot\nabla_{e_j}h(e_k)\cdot\psi+\sum_{j=1}^n\xi\cdot e_j\cdot\nabla_\xi h(e_j)\cdot\psi
\ee
and for every $Z\in \Gamma(Q)$,
\be
\mathrm{Ric}_M(Z)\cdot\psi&=&-4\alpha(h(Z)+\beta Z)\cdot\xi\cdot\psi+2h^2(Z)\cdot\psi+4(n-1)\beta^2 Z\cdot\psi\\
& &+\frac{1}{2}\sum_{j=1}^n\xi\cdot e_j\cdot\nabla_Zh(e_j)\cdot\psi-\sum_{j=1}^n\xi\cdot e_j\cdot\nabla_{e_j}h(Z)\cdot\psi\\
& &-\nabla_\xi h(Z)\cdot\psi+g(Z,\kappa)(-2\alpha\psi+\xi\cdot\Omega\cdot\psi-\kappa\cdot\psi)\\
& &+\nabla_Z^M\kappa\cdot\psi-h(Z)\cdot\xi\cdot\kappa\cdot\psi.
\ee
Furthermore, one has
\be
\mathrm{Scal}_M\psi&=&(4n(n-1)\beta^2-|h|^2-2|\kappa|^2)\psi-8n\alpha\beta\xi\cdot\psi+8\alpha\xi\cdot\Omega\cdot\psi\\
& &-\sum_{j,k=1}^n\xi\cdot e_j\cdot e_k\cdot\nabla_{e_k}h(e_j)\cdot\psi+2\sum_{j=1}^ne_j\cdot\nabla_\xi h(e_j)\cdot\psi\\
& &+4\alpha\kappa\cdot\psi+2\xi\cdot\kappa\cdot\Omega\cdot\psi-2\sum_{j=1}^n e_j\cdot\nabla_{e_j}^M\kappa\cdot\psi.
\ee
\end{ethm}

\noindent{\bf Proof.}
Plugging Equations (\ref{eqkilling}) in (\ref{eqrestrictspin}) gives with the use of (\ref{eqGaussQ}) that for $Z,W\in \Gamma(Q)$  
\be
\nabla_Z^M\nabla_W^M\psi&=&\nabla_Z^M\Big(\beta\xi\cdot W\cdot\psi+\frac{1}{2}\xi\cdot h(W)\cdot\psi\Big)\\
&=&\beta\xi\cdot \nabla_ZW\cdot\psi+\frac{1}{2}\xi\cdot h(\nabla_ZW)\cdot\psi+\frac{1}{2}\xi\cdot \nabla_Zh(W)\cdot\psi\\
& &+\frac{\beta}{2}\Big(h(Z)\cdot W\cdot\psi+h(W)\cdot Z\cdot\psi\Big)\\
& &+\frac{1}{4}h(Z)\cdot h(W)\cdot\psi+\beta^2W\cdot Z\cdot\psi.
\ee
By the fact that the torsion of $\nabla^M$ is zero we get from (\ref{eqGaussQ}) that $[Z,W]=\nabla_ZW-\nabla_WZ-2g(h(Z),W)\xi$ so that
\be
\nabla_{[Z,W]}^M\psi
&=&\beta\xi\cdot(\nabla_ZW-\nabla_WZ)\cdot\psi+\frac{1}{2}\xi\cdot h(\nabla_ZW-\nabla_WZ)\cdot\psi\\
& &-2g(h(Z),W)\Big(\alpha\xi\cdot\psi+\frac{1}{2}\Omega\cdot\psi+\frac{1}{2}\xi\cdot\kappa\cdot\psi\Big).
\ee
We deduce that
\be
R_{Z,W}^M\psi&=&\nabla_{[Z,w]}^M\psi-[\nabla_Z^M,\nabla_W^M]\psi\\
&=&-2g(h(Z),w)\Big(\alpha\xi\cdot\psi+\frac{1}{2}\Omega\cdot\psi+\frac{1}{2}\xi\cdot\kappa\cdot\psi\Big)\\
& &-\frac{1}{2}\xi\cdot (\nabla_Zh(W)-\nabla_Wh(Z))\cdot\psi\\
& &+\frac{1}{4}\pa{h(W)\cdot h(Z)\cdot\psi-h(Z)\cdot h(W)\cdot\psi}+\beta^2\pa{Z\cdot W\cdot\psi-W\cdot Z\cdot\psi}.
\ee
On the other hand since $[\Omega,Z]=2h(Z)$, one has
\be
\nabla_Z^M\nabla_\xi^M\psi&=&\nabla_Z^M\Big(\alpha\xi\cdot\psi+\frac{1}{2}\Omega\cdot\psi+\frac{1}{2}\xi\cdot\kappa\cdot\psi\Big)\\
&=&\frac{\alpha}{2}h(Z)\cdot\psi-\alpha\beta Z\cdot\psi+\frac{1}{2}\nabla_Z^M\Omega\cdot\psi+\frac{\beta}{2}\xi\cdot Z\cdot\Omega\cdot\psi\\
& &+\beta\xi\cdot h(Z)\cdot\psi+\frac{1}{4}\xi\cdot h(Z)\cdot\Omega\cdot\psi+\frac{1}{2}\xi\cdot h^2(Z)\cdot\psi\\
& &+\frac{1}{4}h(Z)\cdot\kappa\cdot\psi-\frac{1}{2}g(h(Z),\kappa)\psi+\frac{\beta}{2}\kappa\cdot Z\cdot\psi+\frac{1}{2}\xi\cdot\nabla_Z^M\kappa\cdot\psi
\ee
and, still using (\ref{eqGaussQ}),
\be
\nabla_\xi^M\nabla_Z^M\psi&=&\nabla_\xi^M\Big(\beta\xi\cdot Z\cdot\psi+\frac{1}{2}\xi\cdot h(Z)\cdot\psi\Big)\\
&=&\beta\xi\cdot \nabla_\xi Z\cdot\psi+\frac{1}{2}\xi\cdot h(\nabla_\xi Z)\cdot\psi+\frac{\beta}{2}\kappa\cdot Z\cdot\psi\\
& &+\beta\xi\cdot h(Z)\cdot\psi+\alpha\beta Z\cdot\psi+\frac{\beta}{2}\xi\cdot Z\cdot\Omega\cdot\psi\\
& &+\frac{1}{4}\kappa\cdot h(Z)\cdot\psi+\frac{1}{2}\xi\cdot\nabla_\xi h(Z)\cdot\psi+\frac{1}{2}\xi\cdot h^2(Z)\cdot\psi\\
& &+\frac{\alpha}{2}h(Z)\cdot\psi+\frac{1}{4}\xi\cdot h(Z)\cdot\Omega\cdot\psi.
\ee
Moreover by the vanishing of the torsion on $M$ we get $[Z,\xi]=-\nabla_\xi Z+g(Z,\kappa)\xi$.
We deduce that
\be
R_{Z,\xi}^M\psi&=&-\nabla_{\nabla_\xi Z}^M\psi+g(Z,\kappa)\nabla_\xi^M\psi-[\nabla_Z^M,\nabla_\xi^M]\psi\\
&=&g(Z,\kappa)\pa{\alpha\xi\cdot\psi+\frac{1}{2}\Omega\cdot\psi+\frac{1}{2}\xi\cdot\kappa\cdot\psi}+2\alpha\beta Z\cdot\psi-\frac{1}{2}\nabla_Z^M\Omega\cdot\psi\\
& &-\frac{1}{2}h(Z)\cdot\kappa\cdot\psi-\frac{1}{2}\xi\cdot\nabla_Z^M\kappa\cdot\psi+\frac{1}{2}\xi\cdot\nabla_\xi h(Z)\cdot\psi.
\ee
Recalling that from its def\mbox{}inition $\Omega$ satisf\mbox{}ies $\Omega(\xi,\cdot)=0$ one can compute the Clif\mbox{}ford action of $\nabla^M\Omega$ and f\mbox{}ind that for every $Z\in \Gamma(Q)$,
\[\nabla_Z^M\Omega\cdot\psi=h^2(Z)\cdot\xi\cdot\psi+\frac{1}{2}\sum_{j=1}^ne_j\cdot\nabla_Zh(e_j)\cdot\psi.\] 
We can hence rewrite
\be
R_{Z,\xi}^M\psi&=&g(Z,\kappa)\pa{\alpha\xi\cdot\psi+\frac{1}{2}\Omega\cdot\psi+\frac{1}{2}\xi\cdot\kappa\cdot\psi}+2\alpha\beta Z\cdot\psi-\frac{1}{2}h^2(Z)\cdot\xi\cdot\psi\\
& &-\frac{1}{2}h(Z)\cdot\kappa\cdot\psi-\frac{1}{2}\xi\cdot\nabla_Z^M\kappa\cdot\psi+\frac{1}{2}\xi\cdot\nabla_\xi h(Z)\cdot\psi\\
& &-\frac{1}{4}\sum_{j=1}^ne_j\cdot\nabla_Zh(e_j)\cdot\psi.
\ee
Applying \cite[p.156]{Hij98} for the local orhonormal frame $\{e_j\}_{1 \leq j\leq n+1}$ of $TM$ with $e_{n+1}=\xi$	we obtain 
\be
\mathrm{Ric}_M(Z)\cdot\psi&=&2\sum_{j=1}^{n+1}e_j\cdot R_{Z,e_j}^M\psi\\
&=&-4\alpha(h(Z)+\beta Z)\cdot\xi\cdot\psi+2h^2(Z)\cdot\psi-h(Z)\cdot\xi\cdot\kappa\cdot\psi\\
& &+\frac{1}{2}\sum_{k=1}^n\xi\cdot e_k\cdot\nabla_Zh(e_k)\cdot\psi-\sum_{k=1}^n\xi\cdot e_k\cdot\nabla_{e_k}h(Z)\cdot\psi-\nabla_\xi h(Z)\cdot\psi\\
& &+4(n-1)\beta^2 Z\cdot\psi+g(Z,\kappa)(-2\alpha\psi+\xi\cdot\Omega\cdot\psi-\kappa\cdot\psi)+\nabla_Z^M\kappa\cdot\psi.
\ee
This shows the second identity of Theorem \ref{tcnexistpsi}.
On the other hand
\be
\mathrm{Ric}_M(\xi)\cdot\psi&=&2\alpha\xi\cdot\kappa\cdot\psi+\kappa\cdot\Omega\cdot\psi+4h(\kappa)\cdot\psi+(|h|^2-|\kappa|^2)\xi\cdot\psi+4n\alpha\beta\psi\\
& &+\frac{1}{2}\sum_{j,k=1}^ne_j\cdot e_k\cdot\nabla_{e_j}h(e_k)\cdot\psi-\sum_{j=1}^n\xi\cdot e_j\cdot\nabla_{e_j}^M\kappa\cdot\psi\\
& &+\sum_{j=1}^n\xi\cdot e_j\cdot\nabla_\xi h(e_j)\cdot\psi.
\ee
This shows the f\mbox{}irst identity of Theorem \ref{tcnexistpsi}.
We compute now the action of the scalar curvature of $M$:
\be
\mathrm{Scal}_M\psi&=&-\sum_{j=1}^{n+1}e_j\cdot\mathrm{Ric}_M(e_j)\cdot\psi\\
&=&(4n(n-1)\beta^2-|h|^2-2|\kappa|^2)\psi-8n\alpha\beta\xi\cdot\psi+8\alpha\xi\cdot\Omega\cdot\psi\\
& &-\sum_{j,k=1}^n\xi\cdot e_j\cdot e_k\cdot\nabla_{e_k}h(e_j)\cdot\psi+2\sum_{j=1}^ne_j\cdot\nabla_\xi h(e_j)\cdot\psi\\
& &+4\alpha\kappa\cdot\psi+2\xi\cdot\kappa\cdot\Omega\cdot\psi-2\sum_{j=1}^n e_j\cdot\nabla_{e_j}^M\kappa\cdot\psi.
\ee
This shows the third identity and achieves the proof of Theorem \ref{tcnexistpsi}.
\hfill$\square$

\begin{erem}\label{remprodloc}
{\rm Under the hypotheses of Theorem \ref{tcnexistpsi}, if one furthermore assumes that $\psi$ is a non-zero $(\alpha,\beta)$-transversal Killing spinor with real $\alpha$ and $\beta$, that $\kappa=0$ and $\nabla h=0$ (e.g. if $M$ is a local Riemannian product or if $M$ is Sasakian) then necessarily $\alpha\beta=0$. From Theorem \ref{tcnexistpsi} the scalar curvature of $M$ must indeed satisfy the condition
\[\mathrm{Scal}_M\psi=(4n(n-1)\beta^2-|h|^2)\psi-8n\alpha\beta\xi\cdot\psi+8\alpha\xi\cdot\Omega\cdot\psi,\]
where $\psi$ is the $(\alpha,\beta)$-transversal Killing spinor on $M$. Taking the Hermitian product with $\xi\cdot\psi$ and identifying the real parts one obtains $0=-8n\alpha\beta|\psi|^2$.
Since $\psi$ does not vanish identically we deduce that $\alpha\beta=0$.
}
\end{erem}

\subsection{Examples of transversal Killing spinors} 

\noindent We construct a f\mbox{}irst important family of examples of Riemannian f\mbox{}lows with transversal Killing spinors. Recall for the next proposition that the unit circle $\s^1$ carries two dif\mbox{}ferent spin structures, the trivial one that we call $(\delta=0)$-spin structure and the non-trivial one that we call $(\delta=1)$-spin structure. We also recall that a $\beta$-Killing spinor on a Riemannian spin manifold $N^n$ is a section $\psi$ of $\Sigma N$ satisfying
\[\nabla_X^N\psi=\beta X\buil{\cdot}{N}\psi\]
for every $X\in \Gamma(TN)$. If a non-zero such spinor f\mbox{}ield exists, then $N^n$ is Einstein with scalar curvature $4n(n-1)\beta^2$ (see e.g. \cite[Prop. 5.12]{Hij98} or \cite{BFGK}), hence $\beta$ must be either real or purely imaginary. The classif\mbox{}ication of the Riemannian spin manifolds with non-trivial $\beta$-Killing spinors was achieved in \cite{Wang89, Baer90, Baum892}.

\begin{prop}\label{pS1fibres}
Let $N$ be an $n$-dimensional Riemannian spin manifold carrying a $\beta$-Killing spinor $\psi$ for some $\beta\in\mathbb{C}$ and $M\bui{\lra}{\pi}N$ be a Riemannian submersion which is either a $\s^1$-bundle with totally geodesic f\mbox{}ibres over $N$ or the second projection of the Riemannian product $M:=\R\times N$ onto $N$. Let $M$ be endowed with the spin structure induced by that of $N$ and the trivial spin structure on $\s^1$ or $\R$ respectively. Then the following holds:
\beit\item[i)] The spinor $\psi$ on $N$ induces a $(0,\beta)$-transversal Killing spinor on $M$.
\item[ii)] In the second case ($M=\R\times N$) if moreover $\beta=0$ then $\psi$ induces an $(\alpha,0)$-transversal Killing spinor on $M$ for any $\alpha\in\mathbb{C}$. Furthermore for $\alpha\in\R$ the spinor field $\psi$ descends to the Riemannian product $\s^1\times N$ if and only if $\alpha\in\frac{\pi\delta}{L}+\frac{2\pi\mathbb{Z}}{L}$, where $L$ is the length of the unit circle. 
\eeit
\end{prop}

\noindent{\bf Proof.} We recall the following lemma about spinors on submersions and $\s^1$-bundles, see \cite[Chap. 1]{Moroithese} or \cite{AmmBaer98}:

\begin{elemme}\label{lstructspinS1}
Let $M\bui{\lra}{\pi}N$ be as in {\rm Proposition \ref{pS1fibres}}.
Then the following statements hold:
\ben\item The manifold $M$ defines a minimal spin Riemannian f\mbox{}low w.r.t. the unit fundamental vector f\mbox{}ield $\xi$ given by the $\s^1$-action or $\frac{\partial}{\partial t}$ respectively and carries a spin structure which is induced by those of $N$ and the trivial one on $\s^1$ or $\R$ respectively.
\item The spinor bundle of $Q$ can be identif\mbox{}ied with $\pi^*\Sigma N$, in particular $\xi\cdot X^*\cdot\varphi=X\buil{\cdot}{N}\varphi$ for every $X\in \Gamma(TN)$, where $X^*\in \Gamma(Q)$ denotes the horizontal lift of $X$ to $M$.
\item For every $\varphi\in\Gamma(\Sigma N)$ (which is identif\mbox{}ied to $\varphi\circ\pi\in\Gamma(\pi^*\Sigma N)$) one has
\[\left|\begin{array}{ll}\nabla_{X^*}\varphi&=\nabla_X^N\varphi\\
\nabla_\xi\varphi&=0\end{array}\right. \]
for every $X\in \Gamma(TN)$. Besides a spinor $\phi$ on $M$ is projectable on $N$ if and only if $\nabla_\xi\phi=0$.
\een
\end{elemme}

\noindent {\bf Proof of {\rm Proposition \ref{pS1fibres}} (continued).} Since $\psi$ is a $\beta$-Killing spinor on the base manifold $N$, then we deduce from Lemma \ref{lstructspinS1} that it satisf\mbox{}ies
\[\left|\begin{array}{ll}\nabla_{Z}\psi&=\nabla_Z^N\psi=\beta Z\buil{\cdot}{N}\psi=\beta\xi\cdot Z\cdot\psi\\
\nabla_\xi\psi&=0\end{array}\right. \]
for every $Z\in \Gamma(\pi^*(TN))\cong Q$, hence $\psi$ is a $(0,\beta)$-transversal Killing spinor on $M$. Note that in the case $n$ odd we identify $\psi$ as a section of the f\mbox{}irst component $\Sigma Q$ of $\Sigma M$. This shows $i)$. Assume now that $M:=\R\times N$ and $\beta=0$. For any $\alpha\in\mathbb{C}$ we set 
\[\phi:=\left|\begin{array}{ll}e^{-i\alpha t}\psi_+ +e^{i\alpha t}\psi_-&\textrm{if }n\textrm{ is even}\\
e^{-i\alpha t}(\psi\oplus\psi)&\textrm{if }n\textrm{ is odd}\end{array}\right.\]
where, if $n$ is even, $\psi=\psi_++\psi_-$ is the decomposition of $\psi$ w.r.t. the Clif\mbox{}ford action of $i\xi$, see above. We check that, under the hypotheses of Proposition \ref{pS1fibres}, the spinor $\phi$ is an $(\alpha,0)$-transversal Killing spinor on $M$. We just describe the case $n$ even, the case $n$ odd being completely analogous. Since $\psi_+$ and $\psi_-$ are parallel, as a consequence of $\psi$ being parallel, then  for all $Z\in \Gamma(Q)$ we deduce $\nabla_Z\phi=0$ which is $\beta\xi\cdot Z\cdot\phi$. 
Moreover
\be
\nabla_\xi\phi&=&i\alpha(-e^{-i\alpha t}\psi_+ +e^{i\alpha t}\psi_-)\\
&=&i\alpha(-i\xi\cdot)(e^{-i\alpha t}\psi_+ +e^{i\alpha t}\psi_-)\\
&=&\alpha\xi\cdot\phi,
\ee
hence $\phi$ is an $(\alpha,0)$-transversal Killing spinor. Finally, let $M:=\s^1\times N$, where $\s^1$ carries the (left-, right- or bi-)invariant metric for which $\mathrm{Length}(\s^1)=L>0$ and the $\delta$-spin structure (where $\delta\in\{0,1\}$). If $\alpha$ is real, then the spinor $\phi$ constructed above on $\R\times N$ satisf\mbox{}ies the equivariance condition (see just after Notes \ref{rspinKTrevet} below) $\phi_{(t,x)}=e^{i\pi\delta}\phi_{(t-L,x)}$
for every $(t,x)\in \R\times N$ if and only if $e^{i\pi\delta}e^{-i\alpha L}=1$. In other words, $\phi$ descends to an $(\alpha,0)$-transversal Killing spinor on $\s^1\times N$ for the metric and spin structure above if and only if $\alpha\in\frac{\pi\delta}{L}+\frac{2\pi\mathbb{Z}}{L}$.
This shows $ii)$ and achieves the proof of Proposition \ref{pS1fibres}. 
\hfill$\square$

\begin{erems}\label{rspinKTrevet}
{\rm\noindent\ben\item It follows from Note \ref{remprodloc} that, on a (local) Riemannian product, the exis\-ten\-ce of a non-zero $(\alpha,\beta)$-transversal Killing spinor with real $\alpha$ and $\beta$ implies $\alpha\beta=0$. Hence the hypothesis $\alpha\beta=0$ in Proposition \ref{pS1fibres} cannot be removed. Moreover the restriction on $M$ being the Riemannian pro\-duct $\s^1\times N$ in Proposition \ref{pS1fibres} $ii)$ turns out to be necessary as well: in general there do not exist non-zero $(\alpha,0)$-transversal Killing spinors with non-zero real $\alpha$ on a given $\s^1$-bundle with totally geodesic f\mbox{}ibres over a spin manifold with parallel spinors. For example, Heisenberg manifolds only admit transversally parallel spinors, see Examples \ref{exS1fibrestg} and Note \ref{exHeis} below.
\item Let $M:=\s^1\times N$ as in Proposition \ref{pS1fibres}. Because of Lemma \ref{lstructspinS1} Proposition \ref{pS1fibres} can be applied to the existence of $(0,\beta)$-transversal Killing spinors on $M$ only if $\s^1$ carries the trivial spin structure. Actually if one f\mbox{}ixes the non-trivial spin structure on $\s^1$ (cor\-res\-pon\-ding to $\delta=1$) then $M$ does not carry any non-trivial $(0,\beta)$-transversal Killing spinor. In the case $\beta=0$ this can be already read of\mbox{}f Proposition \ref{pS1fibres} $ii)$ since $0\notin\frac{\pi}{L}+\frac{2\pi\mathbb{Z}}{L}$. However we give a more general argument which works for any $\beta$. Assume the existence of such a spinor f\mbox{}ield $\psi$ on $M$. Fixing a suf\mbox{}f\mbox{}iciently small nonempty open subset $U$ of $N$ one could write on $\s^1\times U$ the spinor $\psi$ as $\psi=\sum_j f_j\psi_j$, where the $\psi_j$'s are local trivializations of $\Sigma N$ and $f_j$ are sections of $\pi_1^*(\Sigma \s^1)$ (here $\pi_1:M\rightarrow \s^1$ denotes the projection onto the f\mbox{}irst factor). Now a section of $\pi_1^*(\Sigma \s^1)$ w.r.t. the $\delta$-spin structure on $\s^1$ should be a smooth map $f:\R\times N\rightarrow\C$ such that $f(t+1,\cdot)=(-1)^\delta f(t,\cdot)$ for every $t\in\R$. Since the $\psi_j$'s do not depend on the f\mbox{}irst factor $\s^1$ one should therefore have $0=\nabla_\xi\psi=\sum_j\xi(f_j)\psi_j$,
hence together with $\delta=1$ we obtain $f_j=0$ for every $j$, which is a contradiction.
\item The only real line bundles $M\bui{\lra}{\pi}N$ which are minimal Riemannian f\mbox{}lows over some Riemannian manifold $N$ are Riemannian products $\R\times N$. Indeed the vector bundle $M\rightarrow N$ should be trivial since it should possess a global nowhere-vanishing smooth section $\xi$. Moreover, the only metric making such a product into a minimal Riemannian f\mbox{}low is the Riemannian product. Therefore we don't restrict the generality when considering Riemannian products $\R\times N$ instead of arbitrary line bundles over $N$.
\een
}
\end{erems}

\begin{eexemples}\label{exS1fibrestg}
{\rm\noindent\ben
\item The most simple examples that come in mind as application of Proposition \ref{pS1fibres} are the Euclidean space and any f\mbox{}lat torus. Since they admit w.r.t. their trivial spin structure parallel spinors, applying Proposition \ref{pS1fibres} one deduces that they carry an $(\alpha,0)$-transversal Killing spinor for a suitable choice of $\alpha$.
\item In the same way any Riemannian product $\R\times \s^n$ or $\s^1\times \s^n$ (for $n>1$), where $\s^n$ carries its canonical metric and $\s^1$ with its trivial spin structure, admits $(0,\pm\frac{1}{2})$-transversal Killing spinors.
\item Any $\s^1$-bundle which is also a Riemannian submersion with totally geodesic f\mbox{}ibres over $\s^n$ carrying its canonical metric of sectional curvature $1$ admits $(0,\pm\frac{1}{2})$-transversal Killing spinors.
\item As a particular case of the preceding example consider the Hopf-f\mbox{}ibration $\s^3\bui{\lra}{\pi}\C\mathrm{P}^1$, where $\s^3$ is the $3$-dimensional Euclidean sphere and $\C\mathrm{P}^1$ is the complex projective space with its Fubini-Study metric of constant holomorphic sectional curvature $4$. It is well-known that $\pi$ is a $\s^1$-bundle and a Riemannian submersion with totally geodesic f\mbox{}ibres. Furthermore the spin structure on $\s^3$ induced by those of $\C\mathrm{P}^1$ and the trivial spin structure on $\s^1$ is its standard spin structure since there is only one spin structure on $\s^3$. Identifying $\C\mathrm{P}^1$ with the $2$-dimensional sphere $\s^2$ together with $\frac{1}{4}\mathrm{can}$ we have a $2$-dimensional space of $\pm 1$-Killing spinors on $\C\mathrm{P}^1$. We deduce from Proposition \ref{pS1fibres} that $\s^3$ carries a $2$-dimensional space of $(0,\pm 1)$-transversal Killing spinors.\\
More generally, every lens space $\mathbb{Z}_k\setminus \s^3$ with its canonical metric and its trivial spin structure is also a $\s^1$-bundle with totally geodesic f\mbox{}ibres over $\C\mathrm{P}^1$, therefore it admits a non-zero $(0,\pm 1)$-transversal Killing spinor.
\item Let 
\[ G:=\left\{\left(\begin{array}{ccc}1& x& z\\ 0 & 1 & y\\ 0& 0& 1\end{array}\right),\qquad x,y,z\in\R\right\} \]
be the \emph{Heisenberg group}, which is a $3$-dimensional non-compact connected non-abelian Lie group. For a f\mbox{}ixed $r\in(\mathbb{Z}\setminus\{0\})$ consider the discrete subgroup $\Gamma_r$ of $G$ def\mbox{}ined by
\[\Gamma_r:=\left\{\left(\begin{array}{ccc}1& rx& z\\ 0 & 1 & y\\ 0& 0& 1\end{array}\right),\qquad x,y,z\in\mathbb{Z}\right\}. \]
The (homogeneous) quotient $M_r:=\vspace{-1mm}{\Gamma_r}\hspace{-1mm}\setminus\hspace{-0.5mm}\vspace{1mm} G$ is a compact $3$-dimensional manifold called a \emph{Heisenberg manifold}. It carries a two-pa\-ra\-me\-ter family of left-invariant Riemannian metrics which make it into a Riemannian $\s^1$-principal bundle with totally geodesic f\mbox{}ibres over a f\mbox{}lat two-dimensional torus $\mathbb{T}^2:=r\mathbb{Z}\oplus\mathbb{Z}\setminus\R^2$ \cite{AmmBaer98}. Fixing a f\mbox{}lat metric and the trivial spin structure on $\mathbb{T}^2$ we have a $2$-dimensional space of parallel spinors on $\mathbb{T}^2$. Hence it follows from Proposition \ref{pS1fibres} that, for the induced metric $g$ and the induced spin structure on $M_r$, there exists a $2$-dimensional space of $(0,0)$-transversal Killing spinors on $(M_r,g)$. This has been already proved by G. Habib in \cite{Habibthese} where the author performs a direct computation.
\item Let $M:=\wit{\mathrm{PSL}_2(\R)}$ be the universal covering of the projective special linear group of $\R^2$. It can be identif\mbox{}ied with the unitary tangent bundle (or, equivalently, the bundle of positively-oriented orthonormal frames) $\mathbb{U}H^2$ of the hyperbolic plane $H^2$. Fixing the canonical metric and spin structure on $H^2$ we have a $2$-dimensional space of $\pm i$-Killing spinors on $H^2$. From Proposition \ref{pS1fibres} we deduce that, for the induced metric and spin structure on $M$ we have a $2$-dimensional space of $(0,\pm i)$-transversal Killing spinors on $M$.
\een}
\end{eexemples}

\section{Transversal Killing spinors on reducible Riemannian flows} \label{exakillspi} 
\subsection{Equivariance conditions on local Riemannian products}

\noindent In this section, we aim to study the existence of transversal Killing spinors on flat or reducible Riemannian flows (local products).\\

\noindent Let $\Gamma$ be a discrete group acting freely, properly and discontinuously on a spin Riemannian manifold $\ovl{M}$. Assume furthermore that its action preserves both the metric and the spin structure of $\ovl{M}$. Then the quotient manifold $\Gamma\setminus \ovl{M}$ inherits from $\ovl{M}$ a metric and a spin structure such that the canonical projection $\ovl{M}\lra \Gamma\setminus \ovl{M}$ is a Riemannian covering map preserving the spin structures. Furthermore spinor f\mbox{}ields on $\Gamma\setminus \ovl{M}$ are in one-to-one correspondence with $\Gamma$-equivariant spinor f\mbox{}ields on $\ovl{M}$.\\\\
In particular if $M=\Gamma\setminus\R^{n+1}$ then denoting by $t:\Gamma\rightarrow\R^{n+1}$ and $r:\Gamma\rightarrow\mathrm{SO}_{n+1}$ the f\mbox{}irst and second projections of $\mathrm{Isom}_+(\R^{n+1})=\R^{n+1}\rtimes\mathrm{SO}_{n+1}$ (group of orientation-preserving isometries of $\R^{n+1}$) respectively, the equivariance conditions above can be rewritten under the form \cite[Prop. 3.2]{Pfaeffle00}
\beit\item The manifold $M$ is spin if and only if there exists a group-homomorphism $\varepsilon:\Gamma\rightarrow\mathrm{Spin}_{n+1}$ such that the following diagram commutes
\[\xymatrix{&\ar[d]^{\mathrm{Ad}}\mathrm{Spin}_{n+1}\\ \Gamma\ar[ur]^{\varepsilon}\ar[r]^{r}& \mathrm{SO}_{n+1}}\]
\item A spinor f\mbox{}ield on $M$ is a smooth map $\ovl{\psi}:\R^{n+1}\rightarrow\C^{2^{[\frac{n+1}{2}]}}$ satisfying
\[\ovl{\psi}_x=\varepsilon(\gamma)\ovl{\psi}_{\gamma^{-1}(x)} \]
for all $x\in\R^{n+1}$ and $\gamma\in\Gamma$.
\eeit
Next we use this equivariance principle to study the f\mbox{}lat or reducible Riemannian f\mbox{}lows carrying transversal Killing spinors. We begin with local Riemannian products:

\begin{prop}\label{pprodriemspinKT}
Let $(M^{n+1},g,\mathcal{F})$ be a spin Riemannian f\mbox{}low given by a unit vector field with $\kappa=0$ and $h=0$, i.e., a local Riemannian product of a $1$-dimensional manifold with an $n$-dimensional one. Assume that $M$ carries a non-zero $(\alpha,\beta)$-transversal Killing spinor for complex numbers $\alpha$ and $\beta$. Then $\alpha\beta=0$ and $M$ is the quotient of some Riemannian product $\R\times\ovl{N}$ where $\ovl{N}$ admits a non-trivial $\beta$-Killing spinor. Moreover $\beta=0$ if and only if $(M,g)$ is Ricci f\mbox{}lat. 
\end{prop}

\noindent{\bf Proof.}
Let $\psi$ be a non-zero $(\alpha,\beta)$-transversal Killing spinor on $M$. From the assumption the f\mbox{}low being a local Riemannian product, i.e. $h=0$ and $\kappa=0$, we obtain applying Theorem \ref{tcnexistpsi} that
\begin{equation}\label{eqcourbprodriemloc}
\left|\begin{array}{lll}
\mathrm{Scal}_M\psi&=&4n(n-1)\beta^2\psi-8n\alpha\beta\xi\cdot\psi\\
\mathrm{Ric}_M(\xi)\cdot\psi&=&4n\alpha\beta\psi\\
\mathrm{Ric}_M(Z)\cdot\psi&=&4(n-1)\beta^2 Z\cdot\psi-4\alpha\beta Z\cdot\xi\cdot\psi
\end{array}\right.
\end{equation}
for every $Z\in \Gamma(Q)$.
Splitting $\psi$ into $\psi=\psi_++\psi_-$ we deduce from the f\mbox{}irst identity of (\ref{eqcourbprodriemloc}) that $\mathrm{Scal}_M\psi_\pm=4n(n-1)\beta^2\psi_\pm\pm 8in\alpha\beta\psi_\pm$ (the identity holds for both $+$ and $-$).
If $\beta\neq 0$ we know from Notes \ref{remdefKT} that $\{x\in M\,/\,(\psi_\pm)_x\neq 0\}$ is dense in $M$, so that $\mathrm{Scal}_M=4n(n-1)\beta^2\pm 8in\alpha\beta$ on $M$, which of course implies $\alpha=0$. This proves $\alpha\beta=0$.\\
Let now $\ovl{M}\lra M$ be the universal covering of $M$ and $\ovl{M}$ carry the induced metric and spin structure. From the hypotheses $h=0$ and $\kappa=0$ we have $\ovl{M}=\R\times \ovl{N}$ with product metric and spin structure, where $\ovl{N}$ is a simply-connected Riemannian spin manifold. Note also that the lift $\ovl{\xi}$ of $\xi$ to $\ovl{M}$ is then given by $\ovl{\xi}=\frac{\partial}{\partial t}$. The spinor $\psi$ lifts to an $(\alpha,\beta)$-transversal Killing spinor $\ovl{\psi}$ on $\ovl{M}$ which is $\pi_1(M)$-equivariant.\\
In the case where $\alpha=0$ we know from Lemma \ref{lstructspinS1} that $\ovl{\psi}$ is projectable, i.e., it comes from a spinor f\mbox{}ield $\varphi$ on $\ovl{N}$. Applying further Lemma \ref{lstructspinS1} as in the proof of Proposition \ref{pS1fibres} one actually shows that $\varphi$ is a $\beta$-Killing spinor on $\ovl{N}$.\\
In the case where $\alpha\neq 0$ necessarily $\beta=0$, which is from (\ref{eqcourbprodriemloc}) equivalent to $M$ being Ricci f\mbox{}lat. We show the existence of a parallel spinor on $\ovl{M}$, or equivalently on $\ovl{N}$. For this we simply use the argument for constructing $(\alpha,0)$-transversal Killing spinors out of transversally parallel ones (see proof  of Proposition \ref{pS1fibres}) in the reverse way and set, for every $(t,x)\in\ovl{M}$:
\[\ovl{\varphi}_{(t,x)}:=e^{i\alpha t}(\ovl{\psi}_+)_{(t,x)}+e^{-i\alpha t}(\ovl{\psi}_-)_{(t,x)},\]
where $i\ovl{\xi}\cdot\ovl{\psi}_\pm=\pm\ovl{\psi}_\pm$. It is a straightforward computation to show that $\ovl{\varphi}$ is transversally parallel on $\ovl{M}$. In particular since its covariant derivative along $\ovl{\xi}$ vanishes it induces a spinor on $\ovl{N}$ which is then parallel from Lemma \ref{lstructspinS1}. This achieves the proof.
\hfill$\square$

\begin{erems}\label{rprodriemloc}
\noindent{\rm\ben\item In the case where $\beta\neq 0$ one can deduce from (\ref{eqcourbprodriemloc}) that the ei\-gen\-spa\-ces of the Ricci tensor of the universal cover $\ovl{M}$ of $M$ are pointwise $\R$ (cor\-res\-pon\-ding to the eigenvalue $0$) and $T\ovl{N}$ (corresponding to the eigenvalue $4(n-1)\beta^2$). Since any isometry of $\ovl{M}$ should preserve the eigenspaces of its Ricci tensor it should preserve the orthogonal splitting $T_{(t,x)}\ovl{M}=\R\oplus T_x\ovl{N}$. From \cite[Lemma 7.1]{Moroi99} any such isometry should thereby be of the form $(\gamma_1,\gamma_2)$ where $\gamma_1$ and $\gamma_2$ are orientation-preserving isometries of $\R$ and $\ovl{N}$ respectively. However the fundamental group $\pi_1(M)$ need not be a product, that is, $M$ need not be isometric to a global product of the form $\R\times N$ or $\s^1\times N$. Consider for example the locally reducible Riemannian f\mbox{}low $M:=\mathbb{Z}\setminus(\R\times \s^3)$ where $\mathbb{Z}$ acts from the left on $\R\times \s^3$ by $n\mapsto(n+\mathrm{Id}_\R,(-1)^n\mathrm{Id}_{\s^3})$. The manifold $M$ is spin and carries exactly two spin structures, which correspond to the two possible lifts of $\mathbb{Z}\rightarrow\mathrm{SO}_4$, $n\mapsto (-1)^n\mathrm{Id}_{\s^3}$ to $\mathrm{Spin}_4=\mathrm{Spin}_3\times\mathrm{Spin}_3$. For each choice of spin structure the space of $\mathbb{Z}$-equivariant $(0,\beta)$-transversal Killing spinors on $\R\times \s^3$ is exactly $2$-dimensional. Nevertheless $M$ is clearly not dif\mbox{}feomorphic to a product.
\item However B. Alexandrov, G. Grantcharov and S. Ivanov showed in \cite{AlexGrantchIv98} that, under the assumptions of Proposition \ref{pprodriemspinKT}, if furthermore $n$ is even, $\alpha=0$, $\beta\neq 0$ and $M$ is compact, then in fact $M$ is dif\mbox{}feomorphic - but not necessarily isometric - to $\s^1\times \s^n$.
\een
}\end{erems}

\noindent In the case of f\mbox{}lat Riemannian f\mbox{}lows we can make more precise statements:

\begin{ecor}\label{cflotsplatsspinKT}
Let $\Gamma$ be a discrete subgroup of orientation-preserving isometries acting freely on $\R^{n+1}$ and $(M:=\Gamma\setminus\R^{n+1},g,\mathcal{F})$ be a flat spin manifold with a minimal Riemannian flow $\mathcal{F}$. Assume that $M$ carries a non-zero $(\alpha,\beta)$-transversal Killing spinor for complex numbers $\alpha$ and $\beta$. Then $h=0$, $\beta=0$ and $\psi$ comes from a smooth spinor f\mbox{}ield $\ovl{\psi}$ on $\R^{n+1}$ satisfying the following equivariance condition:
\beit\item[a)] Case $\alpha=0$: the section $\ovl{\psi}$ is constant on $\R^{n+1}$ and $\ovl{\psi}=\varepsilon(\gamma)\ovl{\psi}$ for every $\gamma\in\Gamma$.
\item[b)] Case $\alpha\neq 0$: there exist two constant sections $\ovl{\psi}_+$ and $\ovl{\psi}_-$ on $\R^{n+1}$ with $\ovl{\psi}=e^{-i\alpha\langle x,\ovl{\xi}\rangle}\ovl{\psi}_++e^{i\alpha\langle x,\ovl{\xi}\rangle}\ovl{\psi}_-$ and $\ovl{\psi}_\pm=e^{\pm i\alpha\langle t(\gamma),\ovl{\xi}\rangle}\varepsilon(\gamma)\ovl{\psi}_\pm$ for every $\gamma\in\Gamma$,
\eeit
where $\varepsilon:\Gamma\lra\mathrm{Spin}_{n+1}$ is the lift of $\Gamma$ giving the spin structure on $M$ and $\ovl{\xi}$ is the lift of $\xi$ to $\R^{n+1}$.
\end{ecor}

\noindent{\bf Proof.}
The universal cover of $M$ is by hypothesis isometric to $\R^{n+1}$, and because the f\mbox{}low is assumed to be minimal the lift $\ovl{\xi}$ of $\xi$ to $\R^{n+1}$ is a Killing vector f\mbox{}ield on $\R^{n+1}$. Since every such f\mbox{}ield on $\R^{n+1}$ should be constant hence parallel we f\mbox{}irst deduce that $h=0$, i.e., $M$ should be a local Riemannian product. Applying Proposition \ref{pprodriemspinKT} to the Ricci f\mbox{}lat f\mbox{}low $M$ we immediatly obtain $\beta=0$.\\
Moreover $\ovl{\xi}$ as well as the lift $\ovl{\psi}$ of $\psi$ to $\R^{n+1}$ should be $\Gamma$-equivariant. For $\ovl{\xi}$ this means $\gamma_*\ovl{\xi}=\ovl{\xi}$ for every $\gamma\in\Gamma$, that is,
\begin{equation}\label{eqequivxi}
r(\gamma)(\ovl{\xi})=\ovl{\xi}.
\end{equation}
For $\ovl{\psi}$ we f\mbox{}irst notice that the equation of $(\alpha,0)$-transversal Killing spinors can be explicitly solved on $\R^{n+1}$. Indeed if one decomposes such a spinor f\mbox{}ield $\ovl{\psi}$ in a basis of constant sections of $\Sigma\R^{n+1}$ one straightforward obtains
\[\ovl{\psi}_x=e^{-i\alpha\langle x,\ovl{\xi}\rangle}\ovl{\psi}_++e^{i\alpha\langle x,\ovl{\xi}\rangle}\ovl{\psi}_-\]
for every $x\in\R^{n+1}$, where $\ovl{\psi}_\pm$ is a constant section with $i\ovl{\xi}\cdot\ovl{\psi}_\pm=\pm\ovl{\psi}_\pm$. The formula holds in particular if $\alpha=0$, in which case $\ovl{\psi}$ is simply a constant section, i.e., parallel.
The equivariance condition on $\ovl{\psi}$ now reads $\ovl{\psi}_x=\varepsilon(\gamma)\ovl{\psi}_{\gamma^{-1}(x)}$. If $\alpha=0$ this is equivalent to $\ovl{\psi}=\varepsilon(\gamma)\ovl{\psi}$
for every $\gamma\in\Gamma$, which proves $a)$.\\
Assume for the rest of the proof $\alpha\neq 0$. Using $\gamma(x)=r(\gamma)(x)+t(\gamma)$, the equivariance condition becomes
\begin{eqnarray}\label{eqequivpsialphaneq0}
\nonumber\ovl{\psi}_x&=&\varepsilon(\gamma)\pa{e^{-i\alpha\langle \gamma^{-1}(x),\ovl{\xi}\rangle}\ovl{\psi}_++e^{i\alpha\langle \gamma^{-1}(x),\ovl{\xi}\rangle}\ovl{\psi}_-}\\
&\bui{=}{(\ref{eqequivxi})}&e^{-i\alpha\langle x-t(\gamma),\ovl{\xi}\rangle}\varepsilon(\gamma)\ovl{\psi}_++e^{i\alpha\langle x-t(\gamma),\ovl{\xi}\rangle}\varepsilon(\gamma)\ovl{\psi}_-.
\end{eqnarray}
For a given $\gamma\in\Gamma$ we claim that $\ovl{\xi}\cdot\varepsilon(\gamma)=\epsilon(\gamma)\varepsilon(\gamma)\cdot\ovl{\xi}$ for some $\epsilon(\gamma)\in\{\pm 1\}$. Consider indeed $\mathrm{Ad}(\ovl{\xi}\cdot\varepsilon(\gamma)\cdot\ovl{\xi}^{-1})\in\mathrm{SO}_{n+1}$. The conjugation by $\ovl{\xi}$ acts on $\R^{n+1}$ through $\mathrm{Id}_{\R\ovl{\xi}}\oplus-\mathrm{Id}_Q$, so that because of (\ref{eqequivxi}) it leaves the eigenspaces of $\mathrm{Ad}(\varepsilon(\gamma))=r(\gamma)$ invariant and therefore it commutes with $r(\gamma)$, hence
\[\mathrm{Ad}(\ovl{\xi}\cdot\varepsilon(\gamma)\cdot\ovl{\xi}^{-1})=r(\gamma)=\mathrm{Ad}(\varepsilon(\gamma))\]
which is the claim. If $\ovl{\xi}\cdot\varepsilon(\gamma)=-\varepsilon(\gamma)\cdot\ovl{\xi}$ for some f\mbox{}ixed $\gamma\in\Gamma$ then the identif\mbox{}ication of the $+$ and $-$ components in (\ref{eqequivpsialphaneq0}) gives
$\ovl{\psi}_\pm=0$ and hence $\ovl{\psi}=0$, contradiction. Therefore we necessarily have $\ovl{\xi}\cdot\varepsilon(\gamma)=\varepsilon(\gamma)\cdot\ovl{\xi}$ for every $\gamma\in\Gamma$. It follows for the identif\mbox{}ication of the $+$ and $-$ components in (\ref{eqequivpsialphaneq0})
\[ e^{-i\alpha\langle x,\ovl{\xi}\rangle}\ovl{\psi}_+=e^{-i\alpha\langle x-t(\gamma),\ovl{\xi}\rangle}\varepsilon(\gamma)\ovl{\psi}_+\quad\textrm{ and }\quad e^{i\alpha\langle x,\ovl{\xi}\rangle}\ovl{\psi}_-=e^{i\alpha\langle x-t(\gamma),\ovl{\xi}\rangle}\varepsilon(\gamma)\ovl{\psi}_-,\]
which shows $b)$ and achieves the proof of Corollary \ref{cflotsplatsspinKT}.
\hfill$\square$

\subsection{Existence of $(\alpha,0)$-transversal Killing spinors in low dimensions} \label{exBieberbach}
In this section, we determine all compact f\mbox{}lat minimal $3$-dimensional Riemannian f\mbox{}lows carrying at least one non-zero $(\alpha,0)$-transversal Killing spinor $\psi$ for a suitable complex number $\alpha$ (remember from Corollary \ref{cflotsplatsspinKT} that $\beta$ should vanish). All such manifolds are of the form $\Gamma\setminus\R^3$ where $\Gamma$ is one of the six Bieberbach groups. In the case $\alpha=0$ the manifold $M$ should carry a non-zero transversally parallel and hence a parallel spinor f\mbox{}ield by Equations (\ref{eqrestrictspin}). Therefore $M$ should be a f\mbox{}lat torus with trivial spin structure \cite[Thm 5.1]{Pfaeffle00}, see also the case $\Gamma=G_1$ below. So we assume $\alpha\neq 0$. We mainly keep the notations of \cite[Thm 2.8 and 3.3]{Pfaeffle00} and for each Bieberbach group $G_i$ we determine the possible f\mbox{}ields $\ovl{\xi}$ on $\R^3$ and express the equivariance condition of Corollary \ref{cflotsplatsspinKT} $b)$.\\

\noindent$\bullet$ {\bf Case of $\Gamma=G_1:$} The group $G_1$ is generated by three translations associated to three linearly independent vectors $a_j$ in $\R^3$. In that case $M$ is a f\mbox{}lat torus and the lift $\varepsilon$ of $G_1$ to $\mathrm{Spin}_3$ is given on the generators by $\varepsilon(a_j):=e^{i\pi\delta_j}$ where $\delta_j\in\{0,1\}$. Since $r(G_1)=\{\mathrm{Id}\}$ the equivariance condition (\ref{eqequivxi}) for $\ovl{\xi}$ is empty, i.e., $\ovl{\xi}$ can be any constant vector of unit length in $\R^3$. On the other hand, we should have from Corollary \ref{cflotsplatsspinKT} $b)$ that $\ovl{\psi}_\pm=e^{\pm i\alpha\langle a_j,\ovl{\xi}\rangle}e^{i\pi\delta_j}\ovl{\psi}_\pm$
for every $j=1,2,3$, hence the existence of a non-trivial solution is equivalent to $e^{\pm i(\alpha\langle a_j,\ovl{\xi}\rangle+\pi\delta_j)}=1$, i.e., to
\begin{equation}\label{eqalphaG1}
 \alpha\langle\ovl{\xi},a_j\rangle\in\pi\delta_j+2\pi\mathbb{Z}\qquad\forall j\in\{1,2,3\}.
\end{equation}
Note that this implies $\alpha\in\R$. We conclude that, for any f\mbox{}ixed basis $\{a_1,a_2,a_3\}$ of $\R^3$ and any spin structure $(\delta_1,\delta_2,\delta_3)\in\{0,1\}^3$ there exists a non-zero $(\alpha,0)$-transversal Killing spinor on $G_1\setminus\R^3$ if and only if the relation (\ref{eqalphaG1}) is satisf\mbox{}ied. In particular the space of $(\alpha,0)$-transversal Killing spinors is of $2$-dimensional.\\

\noindent$\bullet$ {\bf Case of $\Gamma=G_2:$} The group $G_2$ is generated by three translations associated to the three vectors $a_1:=(0,0,H)$, $a_2:=(L,0,0)$, $a_3:=(T,S,0)$ in $\R^3$ and by the orthogonal transformation that we denote by $(A,\frac{a_1}{2})$ which is def\mbox{}ined by $x\mapsto Ax+\frac{a_1}{2}$. Here $A$ is the rotation of angle $\pi$ around the $x_3$-axis and $H,L,T,S$ are real parameters with $H,L,S>0$. In that case the lift $\varepsilon$ of $G_2$ to $\mathrm{Spin}_3$ is given on the generators by $\varepsilon(a_1):=-1$, $\varepsilon(a_2):=e^{i\pi\delta_2}$, $\varepsilon(a_3):=e^{i\pi\delta_3}$ and $\varepsilon((A,\frac{a_1}{2})):=e^{i\pi\delta_1}e_1\cdot e_2\cdot,$ where $\{e_1,e_2,e_3\}$ denotes the canonical basis of $\R^3$ and $(\delta_1,\delta_2,\delta_3)\in\{0,1\}^3$. Since $r(G_2)=\{\mathrm{Id},A\}$ the equivariance condition (\ref{eqequivxi}) for $\ovl{\xi}$ reduces to $A\ovl{\xi}=\ovl{\xi}$, that is $\ovl{\xi}=e_3$ or $-e_3$, so that w.l.o.g. we can f\mbox{}ix $\ovl{\xi}=e_3$. Writing the equivariance condition from Corollary \ref{cflotsplatsspinKT} $b)$ one can show the following:
given $H,L,S,T\in\R$ with $H,L,S>0$ and $(\delta_1,\delta_2,\delta_3)\in\{0,1\}^3$, there exists a non-zero $(\alpha,0)$-transversal Killing spinor on $G_2\setminus\R^3$ if and only if $\delta_2=\delta_3=0$ and $\alpha H\in\pi+2\pi\delta_1+4\pi\mathbb{Z}$. In that case the space of $(\alpha,0)$-transversal Killing spinors is of $2$-dimensional.\\

\noindent$\bullet$ {\bf Case of $\Gamma=G_3:$} The group $G_3$ is generated by three translations associated to the three vectors $a_1:=(0,0,H)$, $a_2:=(L,0,0)$, $a_3:=(-\frac{L}{2},\frac{L\sqrt{3}}{2},0)$ in $\R^3$ and by the orthogonal transformation that we denote by $(A,\frac{a_1}{3})$ which is def\mbox{}ined by $x\mapsto Ax+\frac{a_1}{3}$. Here $A$ is the rotation of angle $\frac{2\pi}{3}$ around the $x_3$-axis and $H,L$ are positive real parameters. In that case the lift $\varepsilon$ of $G_3$ to $\mathrm{Spin}_3$ is given on the generators by $\varepsilon(a_1):=-e^{i\pi\delta_1}$, $\varepsilon(a_2):=1$, $\varepsilon(a_3):=1$ and $\varepsilon((A,\frac{a_1}{3})):=e^{i\pi\delta_1}(\frac{1}{2}+\frac{\sqrt{3}}{2}e_1\cdot e_2\cdot)$, where $\delta_1\in\{0,1\}$. Since $r(G_3)=\{\mathrm{Id},A,A^2\}$ the equivariance condition (\ref{eqequivxi}) for $\ovl{\xi}$ reduces to $A\ovl{\xi}=\ovl{\xi}$, that is $\ovl{\xi}=e_3$ or $-e_3$, and w.l.o.g. we can f\mbox{}ix $\ovl{\xi}=e_3$. Writing the equivariance condition one can show the following:
given $H,L\in\R$ with $H,L>0$ and $\delta_1\in\{0,1\}$, there exists a non-zero $(\alpha,0)$-transversal Killing spinor on $G_3\setminus\R^3$ if and only if $\alpha H\in\pi+3\pi\delta_1+6\pi\mathbb{Z}$. In that case the space of $(\alpha,0)$-transversal Killing spinors is of $2$-dimensional.\\

\noindent$\bullet$ {\bf Case of $\Gamma=G_4:$} The group $G_4$ is generated by three translations associated to the three vectors $a_1:=(0,0,H)$, $a_2:=(L,0,0)$, $a_3:=(0,L,0)$ in $\R^3$ and by the orthogonal transformation that we denote by $(A,\frac{a_1}{4})$ which is def\mbox{}ined by $x\mapsto Ax+\frac{a_1}{4}$. Here $A$ is the rotation of angle $\frac{\pi}{2}$ around the $x_3$-axis and $H,L$ are positive real parameters. In that case the lift $\varepsilon$ of $G_4$ to $\mathrm{Spin}_3$ is given on the generators by $\varepsilon(a_1):=-1$, $\varepsilon(a_2):=e^{i\pi\delta_2}$, $\varepsilon(a_3):=e^{i\pi\delta_2}$ and $\varepsilon((A,\frac{a_1}{4})):=e^{i\pi\delta_1}(\frac{1}{\sqrt{2}}+\frac{1}{\sqrt{2}}e_1\cdot e_2\cdot)$, where $(\delta_1,\delta_2)\in\{0,1\}^2$. Since $r(G_4)=\{\mathrm{Id},A,A^2,A^3\}$ the equivariance condition (\ref{eqequivxi}) for $\ovl{\xi}$ reduces to $A\ovl{\xi}=\ovl{\xi}$, that is $\ovl{\xi}=e_3$ or $-e_3$, and w.l.o.g. we can f\mbox{}ix $\ovl{\xi}=e_3$. Writing the equivariance condition one can show the following: given $H,L\in\R$ with $H,L>0$ and $(\delta_1,\delta_2)\in\{0,1\}^2$, there exists a non-zero $(\alpha,0)$-transversal Killing spinor on $G_4\setminus\R^3$ if and only if $\alpha H\in\pi+4\pi\delta_1+8\pi\mathbb{Z}$. In that case the space of $(\alpha,0)$-transversal Killing spinors is $2$-dimensional.\\

\noindent$\bullet$ {\bf Case of $\Gamma=G_5:$} The group $G_5$ is generated by three translations associated to the three vectors $a_1:=(0,0,H)$, $a_2:=(L,0,0)$, $a_3:=(\frac{L}{2},\frac{L\sqrt{3}}{2},0)$ in $\R^3$ and by the orthogonal transformation that we denote by $(A,\frac{a_1}{6})$ which is def\mbox{}ined by $x\mapsto Ax+\frac{a_1}{6}$. Here $A$ is the rotation of angle $\frac{\pi}{3}$ around the $x_3$-axis and $H,L$ are positive real parameters. In that case the lift $\varepsilon$ of $G_5$ to $\mathrm{Spin}_3$ is given on the generators by $\varepsilon(a_1):=-1$, $\varepsilon(a_2):=1$, $\varepsilon(a_3):=1$ and $\varepsilon((A,\frac{a_1}{6})):=e^{i\pi\delta_1}(\frac{\sqrt{3}}{2}+\frac{1}{2}e_1\cdot e_2\cdot)$, where $\delta_1\in\{0,1\}$. Since $r(G_5)=\{\mathrm{Id},A,A^2,A^3,A^4,A^5\}$ the equivariance condition (\ref{eqequivxi}) for $\ovl{\xi}$ reduces to $A\ovl{\xi}=\ovl{\xi}$, that is $\ovl{\xi}=e_3$ or $-e_3$, and w.l.o.g. we can f\mbox{}ix $\ovl{\xi}=e_3$. Writing the equivariance condition one can show the following: given $H,L\in\R$ with $H,L>0$ and $\delta_1\in\{0,1\}$, there exists a non-zero $(\alpha,0)$-transversal Killing spinor on $G_5\setminus\R^3$ if and only if $\alpha H\in\pi+6\pi\delta_1+12\pi\mathbb{Z}$. In that case the space of $(\alpha,0)$-transversal Killing spinors is $2$-dimensional.\\

\noindent$\bullet$ {\bf Case of $\Gamma=G_6:$} The group $G_6$ is generated by three translations associated to the three vectors $a_1:=(0,0,H)$, $a_2:=(L,0,0)$, $a_3:=(0,S,0)$ in $\R^3$ and by the orthogonal transformations that we denote - in the obvious way, see above - by $(A,\frac{a_1}{2})$, $(B,\frac{a_2+a_3}{2})$ and $(C,\frac{a_1+a_2+a_3}{2})$. Here $A$ (resp. $B$, $C$) is the rotation of angle $\pi$ around the $x_3$-axis (resp. $x_1$-, $x_2$-axis) and $H,L,S$ are positive real parameters. Since $r(G_6)\supset\{A,B,C\}$ the vector f\mbox{}ield $\ovl{\xi}$ should satisfy $A\ovl{\xi}=\ovl{\xi}$, $B\ovl{\xi}=\ovl{\xi}$ and $C\ovl{\xi}=\ovl{\xi}$, therefore it should vanish. This means that $G_6\setminus\R^3$ cannot carry any Riemannian f\mbox{}low, hence this case should be eliminated.\\

\noindent To sum up, each of the Bieberbach manifolds $G_j\setminus\R^3$ for $j=1,\ldots,5$ carries non-trivial $(\alpha,0)$-transversal Killing spinors for a suitable $\alpha\in\R$ and suitable spin structure.
\section{Transversal Killing spinors on Sasakian ma\-ni\-folds} \label{sasaki}

\subsection{Integrability conditions for transversal Killing spinors on Sasakian manifolds}

\noindent Let $(M^{2m+1},g,\xi)$ be a Sasakian manifold, see Def\mbox{}inition \ref{defSasaki}. First note that, if $\psi$ is an $(\alpha,0)$-transversal Killing spinor on a spin Sasakian manifold $(M^{2m+1},g)$, then so is every component $\psi_r$ of $\psi$ under the Clif\mbox{}ford action of $\Omega$ (indeed the Clif\mbox{}ford action of $\xi$ preserves $\psi_r$ and $\nabla\Omega=0$), compare with Notes \ref{remdefKT}.2.\\

\noindent Recall for the following corollary that a Riemannian f\mbox{}low is called \emph{$\eta$-Einstein} \cite{Oku62} if and only if there exist real constants $\lambda,\mu$ on $M$ such that $\mathrm{Ric}_M=\lambda\mathrm{Id}_{TM}+\mu\xi^\flat\otimes\xi$.

\begin{prop}\label{pcnexistpsiSasaki}
Under the hypotheses of {\rm Theorem \ref{tcnexistpsi}}, assume furthermore that $(M^{2m+1},g,\xi)$ is Sasakian, that $\psi\neq 0$ and that $\alpha$ and $\beta$ are real. Then the following holds:
\beit
\item[i)] Either $\alpha=0$ or $\beta=0$. If $\alpha\neq 0$ then either $\psi$ is an eigenspinor for the Clif\mbox{}ford action by $\Omega$ or $m$ is odd and $\psi=\psi_r+\psi_{m-r}$ for some $r\in\{0,\ldots,m\}$. 
\item[ii)] If $\alpha=0$ then $(M^{2m+1},g,\xi)$ is $\eta$-Einstein.
\item[iii)] If $\alpha\neq 0$ then $g(\mathrm{Ric}_M(Z),h(Z))=0$ for every $Z\in \Gamma(Q)$. If furthermore $\psi_0\neq 0$ or $\psi_{m}\neq0$ then $(M^{2m+1},g,\xi)$ is $\eta$-Einstein. This happens in particular if $\mathrm{dim}(M)=3$.
\eeit
\end{prop}

\noindent{\bf Proof.} Since on a Sasakian manifold $\kappa=0$ and $\nabla h=0$, we f\mbox{}irst deduce from Note \ref{remprodloc} that $\alpha\beta=0$. The identities proved in Theorem \ref{tcnexistpsi} then simplify to
\begin{equation}\label{eqcourbSasaki}
\left|\begin{array}{lll}
\mathrm{Ric}_M(Z)\cdot\psi&=&-4\alpha(h(Z)+\beta Z)\cdot\xi\cdot\psi+(4(2m-1)\beta^2-2) Z\cdot\psi\\
\mathrm{Ric}_M(\xi)\cdot\psi&=&2m\xi\cdot\psi\\
\mathrm{Scal}_M\psi&=&2m(4(2m-1)\beta^2-1)\psi+8\alpha\xi\cdot\Omega\cdot\psi,
\end{array}\right.
\end{equation}
for every $Z\in \Gamma(Q)$. On every Sasakian manifold one has $\mathrm{Ric}_M(\xi)=2m\xi$ \cite{BGal99} hence the second equation above is trivial.
Consider now the last equation involving the scalar curvature of $M$. Decompose $\psi=\sum_{r=0}^{m}\psi_r$ according to (\ref{eqdecSigmaSasaki}) one obtains with the use of $\xi\cdot\psi_r=(-1)^{r+1}i\psi_r$ for every $r\in\{0,\ldots,m\}$ that  
\begin{equation}\label{eqscalpsirSasaki}
\mathrm{Scal}_M\psi_r=2m(4(2m-1)\beta^2-1)\psi_r+(-1)^r8\alpha(2r-m)\psi_r.
\end{equation} 
We consider two cases:
\noindent\beit\item If $\alpha=0$ then coming back to the f\mbox{}irst equation in (\ref{eqcourbSasaki}) we obtain
\[\mathrm{Ric}_M(Z)=(4(2m-1)\beta^2-2)Z,\]
for every $Z\in \Gamma(Q)$, hence $(M^{2m+1},g,\xi)$ is $\eta$-Einstein and $ii)$ is proved.
\item If $\alpha\neq 0$ then $\beta=0$ and one obtains from (\ref{eqscalpsirSasaki}) 
\[\mathrm{Scal}_M=-2m+(-1)^r8\alpha(2r-m)\]
for every $r\in\{0,\ldots,m\}$ for which $\psi_r$ does not vanish. If
there is more than one such $r$, say $r'$, then one has in particular
$(-1)^r(2r-m)=(-1)^{r'}(2r'-m)$.
If $r+r'\equiv 0\,(2)$ then $r=r'$, contradiction, hence $r+r'\equiv 1\,(2)$, from which one deduces that $r+r'=m$. Therefore such an $r'$ must then be unique (equal to $m-r$) and $m$ should be odd. We have proved $i)$.\\
As for the Ricci tensor on $Q$ in that case, we have the equation
\begin{equation}\label{eqRicsasaki}
 \mathrm{Ric}_M(Z)\cdot\psi_r=4(-1)^r i\alpha h(Z)\cdot\psi_r-2Z\cdot\psi_r
\end{equation}
for every $Z\in \Gamma(Q)$ and every $r\in\{0,\ldots,m\}$ for which $\psi_r$ does not vanish. Taking the Hermitian product of that equation with $h(Z)\cdot\psi_r$ and identifying the real parts one obtains $g(\mathrm{Ric}_M(Z),h(Z))=0$,
and this holds for every $Z\in \Gamma(Q)$. Now if one furthermore assumes that $r=0$ or $r=m$, then  with the use of (\ref{eqidentsasa}) 
one deduces from (\ref{eqRicsasaki}) that, in the case $\psi_0\neq 0,$ that $\mathrm{Ric}_M(Z)=(-2-4\alpha)Z$, and in the case $\psi_{m}\neq 0$, that $\mathrm{Ric}_M(Z)=(-2+4(-1)^{m}\alpha)Z$. Hence the flow is $\eta$-Einstein in that case as well. Note that, if $m$ is odd, then both last expressions of the Ricci tensor are the same, which one could expect since in that case both $\psi_0$ and $\psi_{m}$ could be non-vanishing sections.\\
For $m=1$ the only possible values of $r$ are $0$ and $1$, hence the flow must always be $\eta$-Einstein. This shows $iii)$ and achieves the proof.
\hfill$\square$
\eeit

\begin{erem}\label{exHeis}
{\rm\noindent Consider a Heisenberg manifold $(M_r,g,\xi)$ with metric and spin structure as in Examples \ref{exS1fibrestg}. It is a Sasakian manifold. We have proved in Examples \ref{exS1fibrestg}.5 that $(M_r,g,\xi)$ admits transversally parallel spinors. Actually there is no (non-zero) $(\alpha,\beta)$-transversal Killing spinors on $(M_r,g,\xi)$ for real $(\alpha,\beta)\neq(0,0)$. Assume indeed that $\psi$ were such a spinor f\mbox{}ield. If $\beta\neq 0$ then $\alpha=0$ and from Lemma \ref{lstructspinS1} the spinor field $\psi$ would descend to a $\beta$-Killing spinor on a f\mbox{}lat two-torus with trivial spin structure, contradiction. If $\alpha\neq 0$ then using $\mathrm{Ric}_M=-2\mathrm{Id}_{TM}+4\xi^\flat\otimes\xi$ on $(M_r,g,\xi)$ one would straightforward deduce from (\ref{eqcourbSasaki}) that $\alpha h(Z)\cdot\xi\cdot\psi=0$ for every $Z\in \Gamma(Q)$, contradiction. Therefore $\alpha$ and $\beta$ necessarily vanish.
}
\end{erem}

\subsection{Killing vs. transversal Killing spinors}
We now establish a relation between transversal and ``classical'' Killing spinors on Sasakian manifolds. Recall that a $\mathcal{D}$-homothetic deformation of a given metric $g$ on a Riemannian flow $(M,g,\mathcal{F})$ is a metric of the form $\ovl{g}_t:=t^2g_{|_{\R\xi}}+tg_{|_{Q}}$ for some real number $t>0$.
 
\begin{prop}\label{pSasakiKilling}
Let $(M^{2m+1},g,\xi)$ be a spin Sasakian manifold.
\beit\item[a)] The space of $-\frac{1}{2}$-Killing spinors in $\Sigma_0M$ coincides with that of $(-\frac{m+1}{2},0)$-transversal Killing spinors in $\Sigma_0M.$ In particular, a section $\psi_0$ of $\Sigma_0M$ is an $(\alpha,0)$-transversal Killing spinor on $(M^{2m+1},g,\xi)$ for some $\alpha<0$ if and only if there exists a $\mathcal{D}$-homothetic deformation $\ovl{g}_t$ of $g$ for which the corresponding  spinor f\mbox{}ield $\ovl{\psi_0}^t$ is a $-\frac{1}{2}$-Killing spinor.
\item[b)] The space of $\frac{(-1)^m}{2}$-Killing spinors in $\Sigma_mM$ coincides with that of \break $((-1)^m \frac{m+1}{2},0)$-transversal Killing spinors in $\Sigma_m M$. In particular, a section $\psi_m$ of $\Sigma_mM$ is an $(\alpha,0)$-transversal Killing spinor on $(M^{2m+1},g,\xi)$ for some $\alpha$ such that $(-1)^m\alpha >0$ if and only if there exists a $\mathcal{D}$-homothetic deformation $\ovl{g}_t$ of $g$ for which the corresponding spinor f\mbox{}ield $\ovl{\psi_m}^t$ is a $\frac{(-1)^m}{2}$-Killing spinor. 
\eeit
\end{prop}

\noindent{\bf Proof.} First remember that, if $\beta=0$, then every component $\psi_r$ of $\psi$ is again an $(\alpha,0)$-transversal Killing spinor on $M$, therefore we may talk about transversal Killing spinors lying in one of the components $\Sigma_rM$ of $\Sigma M$. Using (\ref{eqrestrictspin}) we compare $\nabla^M\varphi_r$ with $\nabla\varphi_r$ for any section $\varphi_r$ of $\Sigma_rM$: on the one hand
\begin{eqnarray}\label{eqnablaMxipsir}
\nonumber\nabla_\xi^M\varphi_r&=&\nabla_\xi\varphi_r+\frac{1}{2}\Omega\cdot\varphi_r+\frac{1}{2}\xi\cdot\kappa\cdot\varphi_r\\
&=&\nabla_\xi\varphi_r-(-1)^r(r-\frac{m}{2})\xi\cdot\varphi_r,
\end{eqnarray}
and on the other hand, for every $Z\in \Gamma(Q)$,
\begin{equation}\label{eqnablaMXpsir}
\nabla_Z^M\varphi_r=\nabla_Z\varphi_r+\frac{1}{2}\xi\cdot h(Z)\cdot\varphi_r.
\end{equation}
For $r=0$ the identity (\ref{eqnablaMxipsir}) becomes
\[\nabla_\xi^M\varphi_0=\nabla_\xi\varphi_0+\frac{m}{2}\xi\cdot\varphi_0 \]
and for the identity (\ref{eqnablaMXpsir}) we write
\be
\nabla_Z^M\varphi_0&=&\nabla_Z\varphi_0-\frac{1}{2} h(Z)\cdot\xi\cdot\varphi_0\\
&=&\nabla_Z\varphi_0-\frac{1}{2}Z\cdot\varphi_0
\ee
for every $Z\in \Gamma(Q)$. So that bringing together (\ref{eqnablaMxipsir}) and (\ref{eqnablaMXpsir}) we deduce that the spinor f\mbox{}ield $\varphi_0$ is a $-\frac{1}{2}$-Killing spinor on $M$ if and only if it is a $(-\frac{m+1}{2},0)$-transversal Killing spinor. Now if $\alpha<0$ there exists a $t>0$ such that $\frac{\alpha}{t}=-\frac{m+1}{2}$ so that from Lemma \ref{lexistspinKTdeform} there exists a $\mathcal{D}$-homothetic deformation $\ovl{g}_t$ of the metric $g$ for which the corresponding spinor f\mbox{}ield $\ovl{\psi_0}^t$ is a $(-\frac{m+1}{2},0)$-transversal Killing spinor and hence a $-\frac{1}{2}$-Killing spinor on $(M,\ovl{g}_t)$. Furthermore from the argument above the space of $-\frac{1}{2}$-Killing spinors in $\Sigma_0M$ exactly coincides with that of $(-\frac{m+1}{2},0)$-transversal Killing spinors in $\Sigma_0M$. This proves $i)$. For $r=m$ the proof is completely analogous.
\hfill$\square$

\begin{erems}\label{rcnexistpsiSasaki}{
\rm The identities (\ref{eqcourbSasaki}) in the proof of Proposition \ref{pcnexistpsiSasaki} actually provide a link between the sign of $\alpha$ and the geo\-me\-try of $M$. For example if $m$ is odd then the condition $\psi_0\neq 0$ (or alternatively $\psi_m\neq 0$, which gives the same result) implies from (\ref{eqcourbSasaki}) $\mathrm{Scal}_M=-2m-8\alpha m$. 
If $\mathrm{Scal}_M>0$ then necessarily $\alpha<0$, in particular there is no non-trivial such $(\alpha,0)$-transversal Killing spinor for some positive $\alpha$ if $\mathrm{Scal}_M>0$ and $m$ is odd. This will be illustrated with the Berger spheres in the next section.
}
\end{erems}

\begin{ecor}\label{ccsexistpsiSasaki}
Let $(M^{2m+1},g,\xi)$ be a simply-connected complete spin Sasakian manifold carrying a non-zero $(\alpha,\beta)$-transversal Killing spinor $\psi$ for real $\alpha$ and $\beta$. Then we have:
\beit\item[i)] In the case $\alpha=0$ the following holds:
\beit\item If $m=1$ then the manifold $M$ is isometric to $(\s^3,\mathrm{can}),$ up to $\mathcal{D}$-homothetic deformation of $g$, in the case $\beta\neq 0$ and is diffeomorphic to $\R^3$ in the case $\beta=0$.
\item If $m>1$ then $\beta=0$, i.e., $\psi$ is a trans\-ver\-sal\-ly parallel spinor on $M$.
\eeit
\item[ii)] In the case $\alpha\neq 0$ the following holds:
\beit\item If $m=1$ then up to $\mathcal{D}$-homothetic deformation of $g$ the manifold $M$ is isometric to $(\s^3,\mathrm{can})$ if $\alpha<0$ and should satisfy $\mathrm{Ric}_M=-3\mathrm{Id}_{TM}+5\xi^\flat\otimes\xi$ if $\alpha>0$.
\item If $m$ is even, $\psi_0\neq 0$ and $\alpha<0$ (or $\psi_m\neq 0$ and $\alpha>0$ respectively) then $M$ is compact and up to $\mathcal{D}$-homothetic deformation of $g$ it is Einstein-Sasakian.
\item If $m\geq 3$ is odd, $\psi_0+\psi_m\neq 0$ and $\alpha<0$ then $M$ is compact and up to $\mathcal{D}$-homothetic deformation of $g$ it is Einstein-Sasakian or $3$-Sasakian.
\eeit
\eeit
\end{ecor}

\noindent{\bf Proof.} From Proposition \ref{pcnexistpsiSasaki} we know that $\alpha\beta=0$. In the case where $\alpha=0$ and $\beta\neq 0$,  
it follows from the identities (\ref{eqcourbSasaki}) that 
$M$ is Einstein if and only if $\beta^2=\frac{m+1}{2(2m-1)}$. Obviously there exists a $t>0$ such that $\frac{\beta^2}{t}=\frac{m+1}{2(2m-1)}$ so that we deduce the existence of a $\mathcal{D}$-homothetic deformation $\ovl{g}_t$ of $g$ for which $(M,\ovl{g}_t)$ is Einstein with positive Ricci curvature. Since $(M,\ovl{g}_t)$ is moreover complete (for $(M,g)$ is complete if and only if $(M,\ovl{g}_t)$ is) it is compact. One can now adapt an argument {\sl \`a la Hijazi} to show the non-existence of non-zero basic Killing spinors on $(2m+1\geq 5)$-dimensional compact Riemannian f\mbox{}lows with a transversal K\"ahler structure, see e.g. \cite[Thm 5.22]{Hij98}.
Therefore, for $m>1$ then necessarily $\beta=0$. Now 
if $m=1$ and $\beta\neq0$, we deduce with the fact $(M,\ovl{g}_t)$ is of constant curvature, since it is Einstein, that $(M^3,\ovl{g}_t)$ is isometric to $(\s^3,\mathrm{can})$ \cite{Bel00}. 
If $m=1$ and $\beta=0$ then using (\ref{eqcourbSasaki}) the transversal Ricci curvature vanishes and the manifold $M$ is diffeomorphic to $\R^3$ \cite{Blu81}. This shows $i)$.\\
Assume now $\alpha\neq 0$. Then from Proposition \ref{pcnexistpsiSasaki} $\beta=0$ and if we assume moreover that $\psi_0\neq 0$ and $\alpha<0$ (resp. $\psi_m\neq 0$ and $(-1)^m\alpha>0$) then from Proposition \ref{pSasakiKilling} $a)$ there exists a $\mathcal{D}$-homothetic deformation $\ovl{g}_t$ of $g$ for which $\ovl{\psi_0}^t$ is a $-\frac{1}{2}$-Killing spinor on the Sasakian manifold $(M^{2m+1},\ovl{g}_t,\ovl{\xi}^t)$. In particular it is Einstein with scalar curvature $2m(2m+1)$ and since $\ovl{g}_t$ is complete $M$ must actually be compact. 
In the case where $m$ is odd both conditions are equivalent ($\alpha<0$). Furthermore if $m\geq 3$ then it follows from C. B\"ar's classif\mbox{}ication \cite{Baer90} that $(M,\ovl{g}_t)$ should be either Einstein-Sasakian or $3$-Sasakian. If $m=1$ the condition $\psi_0+\psi_1\neq 0$ is fulf\mbox{}illed by hypothesis; if $\alpha<0$ then applying again Proposition \ref{pSasakiKilling} we obtain a $\mathcal{D}$-homothetic deformation $\ovl{g}_t$ of $g$ for which $(M,\ovl{g}_t)$ carries a non-zero $-\frac{1}{2}$-Killing spinor. Hence similarly one concludes that $(M^3,\ovl{g}_t)=(\s^3,\mathrm{can})$. If $m=1$ and $\alpha>0$ then from the identities (\ref{eqcourbSasaki}) the Ricci curvature is given by $\mathrm{Ric}_M=-(4\alpha+2)\mathrm{Id}_{TM}+4(\alpha+1)\xi^\flat\otimes\xi$, therefore $(M^3,\ovl{g}_t)$ satisf\mbox{}ies $\mathrm{Ric}_M=-3\mathrm{Id}_{TM}+5\xi^\flat\otimes\xi$ for some $t>0$. This shows $ii)$ and achieves the proof.
\hfill$\square$

\begin{erem}
{\noindent\rm It also follows from C. B\"ar's classif\mbox{}ication \cite{Baer90} that, conversely, if $(M^{4l+1},g,\xi)$ is a complete simply-connected Einstein-Sasakian manifold, then $M$ is spin and carries non-trivial Killing spinors associated to positive and negative real constants. In the case where $m$ is even we therefore obtain examples of such Riemannian f\mbox{}lows with non-trivial transversal Killing spinors as soon as e.g. one of those constants can be chosen to be $-\frac{1}{2}$ and the corresponding Killing spinor lies pointwise in $\Sigma_0M$. One can obtain such examples when $m\geq 5$ is odd in an analogous way.
}
\end{erem}

\subsection{Example: transversal Killing spinors on the Berger spheres}\label{subsBerger}

\noindent For a positive integer $m$ let $M:=\s^{2m+1}$ be the $(2m+1)$-dimesional sphere equipped with the round metric $g$ with sectional curvature $1$ and its canonical spin structure. It is a Sasakian manifold w.r.t. the vector f\mbox{}ield $\xi_x:=ix$ for every $x\in \s^{2m+1}$, where $\s^{2m+1}\subset\C^{m+1}$ and $i^2=-1$. We shall call the $\mathcal{D}$-homothetic deformations of that Sasakian ma\-ni\-fold the \emph{Berger spheres}. Note that the usual convention is to def\mbox{}ine a Berger metric on $\s^{2m+1}$ as $tg_{|_{\R\xi}}+g_{|_{Q}}$ for some $t>0$.\\
Let the orientation of $\s^{2m+1}$ be such that for every positively oriented basis $\{e_1,\ldots,e_{2m+1}\}$ of $T_x\s^{2m+1}$, the basis $\{x,e_1,\ldots,e_{2m+1}\}$ is positively oriented in $\C^{m+1}$. Choose $\nu_x:=x$ as unit normal on $\s^{2m+1}$. Then one can identify $h$ with $J$ (the standard complex structure on $\C^{m+1}$ restricted to $Q$).

\begin{prop}\label{pspinKTsBerger}
There exists a $1$-dimensional space of $(-\frac{m+1}{2},0)$-(resp. of $((-1)^m\frac{m+1}{2},0)-$) transversal Killing spinors on $\s^{2m+1}$ lying pointwise in 
$\Sigma_0 M$ (resp. $\Sigma_m M$).
\end{prop}

\noindent{\bf Proof.} It is elementary to show that the space of $-\frac{1}{2}$-Killing spinors on $\s^{2m+1}$ lying pointwise in $\Sigma_0M$ is one-dimensional, as well as the space of $\frac{(-1)^m}{2}$-ones in $\Sigma_m M$. The result is then a straightforward consequence of Proposition \ref{pSasakiKilling}. 
\hfill$\square$\\

\noindent There are no other transversal Killing spinors on the Berger spheres as those that have already been constructed: this is the statement of the following proposition, of which proof is left to the reader.

\begin{prop}
Let $m\geq 1$ and assume the existence of a non-zero $(\alpha,\beta)$-transversal Killing spinor $\psi$ on $(\s^{2m+1},g)$ for complex $\alpha,\beta$. Then $\alpha\beta=0$ and
\beit\item[i)] if $\alpha=0$ then $m=1$, $\beta^2=1$ and $\psi$ is one of the spinors constructed in {\rm Examples \ref{exS1fibrestg} 4}. 
\item[ii)] if $\alpha\neq 0$ then $\beta=0$, $\alpha=\epsilon\frac{m+1}{2}$ for some $\epsilon\in\{\pm 1\}$ and $\psi$ is one of the spinors constructed in {\rm Proposition \ref{pspinKTsBerger}}.
\eeit
\end{prop}

\begin{erems}\label{remsphereBerger}
{\rm\noindent\ben\item There exists in particular no non-zero $(\frac{m+1}{2},0)$-transversal Killing spinor on $(\s^{2m+1},g)$ with $m$ odd, although the space of $(-\frac{m+1}{2},0)$-transversal Killing spinors is $2$-dimensional (compare with the case $m$ even). In particular the complex number $\alpha$ cannot be arbitrarily changed into $-\alpha$.
\item Let $M:=\Gamma\setminus \s^3$ where $\Gamma$ is a non-trivial f\mbox{}inite subgroup of $\SU_2$. Remember that, denoting by $\SU(2)\bui{\lra}{\Theta}\SO_3$ be the universal covering map of $\SO_3$, the group $\Gamma$ is conjugated to a subgroup of one of the following subgroups of $\SU_2$ (see e.g. \cite[Rem. p.57]{Ammhabil} or \cite[Thm 2.6.7]{Wolf67}, where certain subgroups of $\SU_2$ are obviously missing): the cyclic group of order $k$ ($k\in\mathbb{N}\setminus\{0\}$) generated by the element $\left(\begin{array}{cc}e^{\frac{2i\pi}{k}}&0\\ 0& e^{-\frac{2i\pi}{k}}\end{array}\right)\in\SU_2$, $D_k^*:=\Theta^{-1}(D_k^+)$ (resp. $T^*:=\Theta^{-1}(T^+)$, $O^*:=\Theta^{-1}(O^+)$ and $I^*:=\Theta^{-1}(I^+)$) where $D_k^+$ (resp. $T^+$, $O^+$ and $I^+$) is the group of orientation-preserving isometries of a regular $k$-gon (resp. tetrahedron, octahedron, and icosahedron). Every such quotient endowed with the metric $g$ induced by the standard metric on $\s^3$ is of course again a Sasakian manifold.
Moreover it is spin and carries a spin structure for which the space of $\frac{1}{2}$-Killing spinors on $(M,g)$ is $2$-dimensional resp. a spin structure for which the space of $-\frac{1}{2}$-Killing spinors on $(M,g)$ is $2$-dimensional see \cite[Cor. 5.2.5]{Ammhabil}. Hence there exists for the latter spin structure a $2$-dimensional space of $(-1,0)$-transversal Killing spinors.
\een}
\end{erems}

\section{Transversal Killing spinors on $3$-dimensional f\mbox{}lows}

\subsection{Integrability conditions for transversal Killing spinors on $3$-dimensional f\mbox{}lows}

\noindent In this section we assume that $(M,g,\mathcal{F})$ is a $3$-dimensional Riemannian f\mbox{}low. We f\mbox{}ix the orientation on $ Q$ induced by those of $M$ and $\xi$ (i.e., a basis $\{Z,W\}$ of $ Q$ is oriented w.r.t. that orientation if and only if $\{\xi,Z,W\}$ is oriented as local basis of $TM$).\\\\
A f\mbox{}irst consequence of the dimension of $M$ being $3$ is the existence of an almost-Hermitian structure $J$ on $ Q$ def\mbox{}ined in a local positively-oriented orthonormal basis $\{e_1,e_2\}$ of $Q$ by the matrix
\[J:=\left(\begin{array}{cc}0&-1\\ 1&0\end{array}\right).\]
It is easy to see that $J$ is well-def\mbox{}ined, i.e., doesn't depend on the choice of local basis of $ Q$ (this follows from the fact that $\mathrm{SO}_2$ is abelian). Furthermore, $J$ is ``K\"ahler'' on $ Q$, that is $\nabla J=0$ 
on $M$. Since $h$ is a skew-symmetric tensor, one may write $h$ as
\[h=bJ \]
for some smooth globally def\mbox{}ined function $b:M\lra\R$. We recall that the complex volume form 
$$\omega_3=-\xi\cdot e_1\cdot e_2$$
acts as the identity on the spinor bundle $\Sigma M$. Hence one may identify the Clif\mbox{}ford action of any $2$-form with that of forms of lower degrees. On the one hand we have 
\[\xi\cdot Z\cdot=J(Z)\cdot\]
for every $Z\in \Gamma(Q)$ and on the other hand we also have
\[Z\cdot W\cdot=g(J(Z),W)\xi\cdot-g(Z,W)\mathrm{Id}_{\Sigma M}\]
for all $Z,W\in \Gamma(Q)$. For example, one can identify the Clif\mbox{}ford action of $\Omega$ through that of $\xi$ by 
$\Omega\cdot\varphi=b\xi\cdot\varphi$ for all $\psi\in \Gamma(\Sigma M)$.


\noindent In the following proposition and henceforth we denote by $db_{|_{ Q}}:=\sum_{k=1}^2e_k(b)e_k$ (orthogonal projection of $\mathrm{grad}(b)$ onto $ Q$).

\begin{prop}\label{pcnexistdim3}
Let $(M^3,g,\mathcal{F})$ be a spin Riemannian f\mbox{}low carrying a non-zero $(\alpha,\beta)$-transversal Killing spinor for real $\alpha$ and $\beta$. Then the following holds:
\beit\item[i)] $\alpha=0$ or $\kappa=0$.
\item[ii)] $d\kappa^\flat(e_1,e_2)=2(\xi(b)-4\alpha\beta)$.
\item[iii)] For every $Z\in \Gamma(Q)$,
\[\left|\begin{array}{ll}
\mathrm{Scal}_M&=2\pa{4\beta^2-b^2-4\alpha b-\mathrm{div}^M(\kappa)}\\\\
\mathrm{Ric}_M(\xi)&=(2b^2-\mathrm{div}^M(\kappa))\xi+J\pa{2b\kappa-db_{|_{ Q}}}\\\\
\mathrm{Ric}_M(Z)&=2(2\beta^2-b^2-2\alpha b)Z+(4\alpha\beta-\xi(b))J(Z)\\
&\phantom{=}+g\pa{J(2b\kappa-db_{|_{ Q}}),Z}\xi+\nabla_Z\kappa-g(Z,\kappa)\kappa.
\end{array}\right.\]
\eeit
\end{prop}

\noindent{\bf Proof.} 
We keep the notation $\{e_j\}_{1\leq j\leq 2}$ for a local orthonormal basis of $Q$ and simplify the terms given in Theorem \ref{tcnexistpsi}. On the one hand
\be
\sum_{j=1}^2 e_j\cdot\nabla_{e_j}^M\kappa\cdot\psi
&=&d\kappa^\flat\cdot\psi+(\mathrm{div}^M(\kappa)-|\kappa|^2)\psi-J(\nabla_\xi\kappa)\cdot\psi+b\kappa\cdot\psi.
\ee
On the other hand using $\xi\lrcorner d\kappa^\flat=\nabla_\xi\kappa$, one can straighforward prove that  
\[d\kappa^\flat\cdot\psi=J(\nabla_\xi\kappa)\cdot\psi+d\kappa^\flat(e_1,e_2)\xi\cdot\psi.\]
We deduce that
\[\sum_{j=1}^2 e_j\cdot\nabla_{e_j}^M\kappa\cdot\psi=(\mathrm{div}^M(\kappa)-|\kappa|^2)\psi+d\kappa^\flat(e_1,e_2)\xi\cdot\psi+b\kappa\cdot\psi.\]
Furthermore, since $h=bJ$ and $J$ is $\nabla$-parallel, one has $\nabla_Xh=X(b)J$ for every $X\in \Gamma(TM)$. In particular, if $Z\in \Gamma(Q)$,
\[\left|\begin{array}{ll}\sum_{k=1}^2\xi\cdot e_k\cdot\nabla_Zh(e_k)\cdot\psi&=-2Z(b)\psi\\\\
\sum_{k=1}^2\xi\cdot e_k\cdot\nabla_{e_k}h(Z)\cdot\psi&=g(J(db_{|_{ Q}}),Z)\xi\cdot\psi-Z(b)\psi.\end{array}\right.\]
Moreover $\sum_{j,k=1}^2 e_j\cdot e_k\cdot\nabla_{e_j}h(e_k)\cdot\psi=-2J(db_{|_{ Q}})\cdot\psi$ and $\sum_{j=1}^2\nabla_{e_j}h(e_j)=J(db_{|_{ Q}})$.
In particular
\[\sum_{j,k=1}^2 e_j\cdot e_k\cdot\nabla_{e_k}h(e_j)\cdot\psi=0.\]
Last we compute $\sum_{j=1}^2e_j\cdot\nabla_\xi h(e_j)\cdot\psi=2\xi(b)\xi\cdot\psi$. Now we begin with the proof of the proposition. From Theorem \ref{tcnexistpsi} we have
\be
\mathrm{Scal}_M\psi
&=&2(4\beta^2-b^2-4\alpha b-\mathrm{div}^M(\kappa))\psi\\
& &+\{4(\xi(b)-4\alpha\beta)-2d\kappa^\flat(e_1,e_2)\}\xi\cdot\psi+4\alpha\kappa\cdot\psi.
\ee
Taking the Hermitian scalar product of this last identity with $\psi$ and identifying the real parts we obtain $\mathrm{Scal}_M=2\{4\beta^2-b^2-4\alpha b-\mathrm{div}^M(\kappa)\}$
and for what remains we deduce that $2(\xi(b)-4\alpha\beta)-d\kappa^\flat(e_1,e_2)=0$ and $\alpha\kappa=0$.
In particular either $\alpha=0$ or $\kappa=0$.
Coming back to the equations involving $\mathrm{Ric}_M$, we have on the one hand
\be
\mathrm{Ric}_M(\xi)\cdot\psi
&=&(2b^2-\mathrm{div}^M(\kappa))\xi\cdot\psi+2bJ(\kappa)\cdot\psi-J(db_{|_{ Q}})\cdot\psi,
\ee
from which we deduce $\mathrm{Ric}_M(\xi)=(2b^2-\mathrm{div}^M(\kappa))\xi+J\pa{2b\kappa-db_{|_{ Q}}}$. On the other hand, it also follows from Theorem \ref{tcnexistpsi} that, for every $Z\in \Gamma(Q)$,
\be
\mathrm{Ric}_M(Z)\cdot\psi
&=&2(2\beta^2-b^2-2\alpha b)Z\cdot\psi+(4\alpha\beta-\xi(b))J(Z)\cdot\psi\\
& &+\{2bg(J(\kappa),Z)-g(J(db_{|_{ Q}}),Z)\}\xi\cdot\psi+\nabla_Z\kappa\cdot\psi\\
& &-g(Z,\kappa)\kappa\cdot\psi,
\ee
from which we deduce that
\be
\mathrm{Ric}_M(Z)&=&2(2\beta^2-b^2-2\alpha b)Z+(4\alpha\beta-\xi(b))J(Z)\\
& &+g\pa{J(2b\kappa-db_{|_{ Q}}),Z}\xi+\nabla_Z\kappa-g(Z,\kappa)\kappa.
\ee
Hence the proof of the proposition is achieved. 
\hfill$\square$\\\\

\noindent It follows from Proposition \ref{pcnexistdim3} that either $\alpha=0$ or $\kappa=0$. Let us examine the last condition.

\begin{prop}\label{pcnexistdim3harm}
Let $(M^3,g,\mathcal{F})$ be a spin Riemannian f\mbox{}low carrying a non-zero $(\alpha,\beta)$-transversal Killing spinor for real $\alpha$ and $\beta$. Assume that $\kappa=0$.
\beit\item[i)] For every $Z\in \Gamma(Q)$ one has
\[\left|\begin{array}{ll}
\mathrm{Scal}_M&=2\pa{4\beta^2-b^2-4\alpha b}\\
\mathrm{Ric}_M(\xi)&=2b^2 \xi-J\pa{db_{|_{ Q}}}\\
\mathrm{Ric}_M(Z)&=2(2\beta^2-b^2-2\alpha b)Z-g\pa{J(db_{|_{ Q}}),Z}\xi.
\end{array}\right.\]
\item[ii)] One has
\[\left|\begin{array}{ll}\nabla_\xi db_{|_{ Q}}&=0\\
\mathrm{div}^M(J(db_{|_{ Q}}))&=-8\alpha\beta(3b+2\alpha).
\end{array}\right.\]
\item[iii)] If furthermore $M$ is compact, then $\alpha\beta=0$ and $\xi(b)=0$.
\eeit
\end{prop}

\noindent{\bf Proof.} If $\kappa=0$ then the equations of Proposition \ref{pcnexistdim3} obviously simplify to the equations in $i)$. It follows from those equations that $M$ is $\eta$-Einstein if and only if $b$ is constant. Moreover, one may write the Ricci tensor of $M$ in the following way:
\be
\mathrm{Ric}_M
&=&2(2\beta^2-b^2-2\alpha b)\mathrm{Id}_{TM}+4(b^2+\alpha b-\beta^2)\xi^\flat\otimes\xi\\
& &-J(db_{|_{ Q}})^\flat\otimes\xi-\xi^\flat\otimes J(db_{|_{ Q}}).
\ee
Since the divergence of $\xi$ vanishes by the fact that $h$ is skew-symmetric. Using the identity $\mathrm{div}^M(X^\flat\otimes Y)=\mathrm{div}^M(X)Y-\nabla_X^MY$ for $X,Y\in \Gamma(TM)$
one obtains, on the one hand, $\mathrm{div}^M(\xi^\flat\otimes\xi)=-\kappa=0$. Similarly, one has
\[\left|\begin{array}{ll}
\mathrm{div}^M(\xi^\flat\otimes J(db_{|_{ Q}}))&=-\nabla_\xi J(db_{|_{ Q}})+bdb_{|_{ Q}}\\\\
\mathrm{div}^M(J(db_{|_{ Q}})^\flat\otimes\xi)&=\mathrm{div}^M(J(db_{|_{ Q}}))\xi+bdb_{|_{ Q}}.\end{array}\right.
\]
Therefore we can compute the divergence of $\mathrm{Ric}_M$ and obtain
\be
\mathrm{div}^M(\mathrm{Ric}_M)
&=&2(b+2\alpha)db_{|_{ Q}}+\nabla_\xi J(db_{|_{ Q}})-(4b\xi(b)+\mathrm{div}^M(J(db_{|_{ Q}})))\xi.
\ee
On the other hand we have $d\mathrm{Scal}_M=-4(b+2\alpha)db_{|_{ Q}}-4(b+2\alpha)\xi(b)\xi$.
The identity $\mathrm{div}^M(\mathrm{Ric}_M)=-\frac{1}{2}d\mathrm{Scal}_M$ implies together with the fact that $J$ is parallel w.r.t. the connection $\nabla,$
\[2(b+2\alpha)db_{|_{ Q}}+J(\nabla_\xi db_{|_{ Q}})-(4b\xi(b)+\mathrm{div}^M(J(db_{|_{ Q}})))\xi=2(b+2\alpha)(db_{|_{ Q}}+\xi(b))\xi,\]
that is,
\begin{equation}\label{eqdivRicdim3kappa0}
\left|\begin{array}{ll}\nabla_\xi db_{|_{ Q}}&=0\\
2(3b+2\alpha)\xi(b)+\mathrm{div}^M(J(db_{|_{ Q}}))&=0.
\end{array}\right.
\end{equation}
Now remember that from Proposition \ref{pcnexistdim3} one has $\xi(b)=4\alpha\beta$, since the mean curvature vanishes. Hence the second equation in (\ref{eqdivRicdim3kappa0}) may be rewritten under the form
$\mathrm{div}^M(J(db_{|_{ Q}}))=-8\alpha\beta(3b+2\alpha)$.
For the rest of the proof assume $M$ to be compact. If $\alpha\beta$ did not vanish, then one would get from Stokes Theorem that 
\begin{equation}\label{eqintbvg}
3\int_M bv_g+2\alpha\mathrm{Vol}(M)=0.
\end{equation}
On the other hand, still following from Stokes Theorem, one would have
\[-8\alpha\beta\int_M b(3b+2\alpha) v_g=\int_M b\mathrm{div}^M(J(db_{|_{ Q}}))v_g=\int_M g\pa{db,J(db_{|_{ Q}})}v_g=0,\]
which would imply
\begin{equation}\label{eqintb2vg}
3\int_Mb^2 v_g+2\alpha\int_M bv_g=0.
\end{equation}
Combining (\ref{eqintb2vg}) with (\ref{eqintbvg}) one would have in particular
\[\int_M b^2 v_g=-\frac{2\alpha}{3} \int_M bv_g=\frac{4\alpha^2}{9}\mathrm{Vol}(M)=\frac{(\int_M bv_g)^2}{\mathrm{Vol}(M)},\]
i.e., $\frac{1}{\mathrm{Vol}(M)}\int_M b^2 v_g=\pa{\frac{1}{\mathrm{Vol}(M)}\int_M bv_g}^2$.
But this is the equality-case in Cauchy-Schwarz inequality, so that $b$ should be constant, which in turn would imply $4\alpha\beta=\xi(b)=0$ contradiction.
Therefore, if $M$ is compact, then $\alpha\beta=0$ in particular $\xi(b)=0$. 
\hfill$\square$

\subsection{Compact $\eta$-Einstein $3$-dimensional minimal f\mbox{}lows with transversal Killing spinors}

\noindent In this section we describe all compact $\eta$-Einstein $3$-dimensional minimal f\mbox{}lows with transversal Killing spinors for real constants $\alpha$ and $\beta$.

\begin{ecor}\label{ccnsetaeinsteindim3}
Let $(M^3,g,\mathcal{F})$ be a Riemannian f\mbox{}low carrying a non-zero $(\alpha,\beta)$-transversal Killing spinor for real $\alpha$ and $\beta$. Assume that $\kappa=0$. Then the following propositions are equivalent:
\beit\item[i)] The Riemannian f\mbox{}low $(M^3,g,\xi)$ is $\eta$-Einstein.
\item[ii)] The function $b$ (which is def\mbox{}ined by $h=bJ$) is constant.
\item[iii)] The manifold $M$ is either a local Riemannian product or a Sasakian manifold up to homothety on the metric.
\eeit
\end{ecor}

\noindent{\bf Proof.} The equivalence of $i)$ with $ii)$ is a direct consequence of Proposition \ref{pcnexistdim3harm} $i)$. As for the equivalence of $ii)$ with $iii)$, one should consider the two cases. The first case is where $b=0$ which gives the vanishing of the O'Neill tensor and with the assumption $\kappa=0$ we deduce that $M$ is locally a Riemannian product. The second case is where $b\neq 0$ which implies that the manifold $(M,b^2g,\frac{1}{b}\xi)$ is Sasakian.
\hfill$\square$



\begin{prop}\label{pclassifKTdim3}
Let $(M^3,g,\mathcal{F})$ be a spin compact Riemannian f\mbox{}low carrying a non-zero $(\alpha,\beta)$-transversal Killing spinor. Assume that the f\mbox{}low is minimal, $\eta$-Einstein and that $\alpha,\beta\in\R$. Then $\alpha\beta=0$ and up to homotheties and $\mathcal{D}$-homothetic deformations of the metric $g$ the manifold $M$ is isometric to one of the following:
\beit\item[i)] If $\beta\neq 0$: $\s^1\times \s^2$, $\s^3$, $\mathbb{Z}_k\setminus \s^3$ for some $k$.
\item[ii)] If $\beta=0$: $\Gamma\setminus \s^3$ for some f\mbox{}inite subgroup $\Gamma\subset\SU_2$ (if $\alpha<0$), $\Gamma\setminus \wit{\mathrm{PSL}_2(\R)}$ for some f\mbox{}inite cocompact subgroup of ${\rm Isom}_+(\wit{\mathrm{PSL}_2(\R)})$ (if $\alpha>0$), a Heisenberg manifold $M_r$ (if $\alpha=0$) or a Bieberbach manifold $\Gamma\setminus\R^3$.
\eeit
\end{prop}

\noindent{\bf Proof.}
We already know from Proposition \ref{pcnexistdim3harm} and Corollary \ref{ccnsetaeinsteindim3} that $\alpha\beta=0$ and that $b$ should be constant. Hence for $b=0$, the manifold $M$ is locally a product of two Riemannian manifolds, whereas for $b\neq 0$ it is a Sasakian manifold up to homothety on the metric. We consider the two cases separately:\\

\noindent$\bullet$ {\bf Case where $b=0$ :} It follows from Proposition \ref{pprodriemspinKT} that the universal cover of $M$ is isometric to $\R\times \ovl{N}$ where $\ovl{N}$ is a simply-connected complete Riemannian surface carrying a $\beta$-Killing spinor. 
In dimension $2$ the only such surfaces are - up to homothety on the metric - $\R^2$ (for $\beta=0$) and $\s^2$ (for $\beta\neq 0$), so that $\widetilde{M}$ is isometric to $\R^3$ or to $\R\times \s^2$ (remember that $\beta$ is assumed to be real).
\beit\item In the subcase $\ovl{N}=\s^2$ we have seen in Notes \ref{rprodriemloc} that the fundamental group of $M$ should be embedded in the product $\mathrm{Isom}_+(\R,\mathrm{can})\times\mathrm{Isom}_+(\ovl{N},g_{\ovl{N}})$ where $\mathrm{Isom}_+$ denotes the group of orientation-preserving isometries of the corresponding Riemannian manifold. Since the only orien\-ta\-tion-pre\-ser\-ving isometry subgroup of $\SO_3=\mathrm{Isom}_+(\s^2)$ acting freely on $\s^2$ is the trivial one, we deduce that $\pi_1(M)$ is a (discrete) subgroup of $\R=\mathrm{Isom}_+(\R,\mathrm{can})$, so that $M$ is either isometric to $\R\times \s^2$ (if $\pi_1(M)=\{\Id\}$) - which is excluded because of $M$ being assumed to be compact - or to $\s^1\times \s^2$ (if $\pi_1(M)\cong\mathbb{Z}$), and in the last situation $\s^1$ carries the trivial spin structure.
\item In the subcase $\ovl{N}=\R^2$, i.e. $\beta=0$, the manifold $M$ is Ricci f\mbox{}lat hence f\mbox{}lat, therefore it is isometric to the quotient $\Gamma\setminus\R^3$ where $\Gamma\subset\mathrm{Isom}_+(\R^3,\mathrm{can})=\R^3\rtimes\SO_3$ is a discrete subgroup of orientation-preserving isometries ac\-ting freely on $\R^3$. In other words, $M$ is one of the Bieberbach manifolds discussed in Example \ref{exBieberbach}.
\eeit
\noindent $\bullet$ {\bf Case where $b\neq 0$}: Up to changing $g$ into $b^2g$ we may assume that $b=1$ so that $M$ is Sasakian. In that case the assertion follows straightforward from Belgun's uniformization theorem \cite{Bel00} stating that $M$ should be a compact quotient of $\s^3$, $\mathrm{Nil}_3$ or $\wit{\mathrm{PSL}_2(\R)}$. This achieves the proof.
\hfill$\square$


\providecommand{\bysame}{\leavevmode\hbox to3em{\hrulefill}\thinspace}





\end{document}